\documentclass[11pt,reqno]{amsart}

\makeatletter

\usepackage{amssymb}
\usepackage{latexsym}
\usepackage{amsbsy}
\usepackage{amsfonts}
\usepackage{hyperref}

\def\marginpar#1{\ignorespaces}

\newtheorem{theorem}{Theorem}[section]
\newtheorem{proposition}[theorem]{Proposition}

\newtheorem{lemma}[theorem]{Lemma}
\newtheorem{corollary}[theorem]{Corollary}
\newtheorem{definition}[theorem]{Definition}

\theoremstyle{definition}
\newtheorem{remark}[theorem]{Remark}

\newtheorem{example}[theorem]{Example}

\numberwithin{equation}{section}

\def\AArm{\fam0 \rm}%
\newdimen\AAdi%
\newbox\AAbo%
\def\AAk#1#2{\setbox\AAbo=\hbox{#2}\AAdi=\wd\AAbo\kern#1\AAdi{}}%

\newcommand{\BBone}{{\ensuremath{{\AArm 1\AAk{-.8}{I}I}}}}

\def\eqref#1{(\ref{#1})}
\def\eqlabel#1{\def\@currentlabel{#1}}

\def\formula#1{\def\@tempa{#1}\let\@tempb\theequation\def\theequation{%
\hbox{#1}}\def\@currentlabel{(\theequation)}$$}
\def\endformula{\leqno\hbox{(\@tempa)}$$\@ignoretrue\let\theequation\@tempb}

\def\given{\hskip5\p@\relax\vrule\@width.4\p@\hskip5\p@\relax}

\newcommand{\open}[1]{%
\par\normalfont\topsep6\p@\@plus6\p@\trivlist\item[\hskip\labelsep\itshape#1%
\@addpunct{.}]\ignorespaces}

\DeclareRobustCommand{\close}[1]{%
  \ifmmode 
  \else \leavevmode\unskip\penalty9999 \hbox{}\nobreak\hfill
  \fi
  \quad\hbox{$#1$}}

\newlength{\toskip}

\def \R {{\mathbb R}}
\def \Q {{\mathbb Q}}

\def \P {{\mathbb P}}
\def \E {{\mathbb E}}

\def \L {{\mathbb L}}

\def \W {{\mathbb W}}

\def \phi {\varphi}

\newcommand     {\NN}{\mathbb{N}}
\newcommand     {\RR}{\mathbb{R}}
\newcommand     {\PP}{\mathbb{P}}
\newcommand     {\EE}{\mathbb{E}}

\makeatother

\begin{document}
\date{\today}

\title[Quasi-stationarity for population diffusion
 processes]{QUASI-STATIONARY
 DISTRIBUTIONS AND DIFFUSION MODELS IN POPULATION DYNAMICS}

 \author[P. Cattiaux]{\textbf{\quad {Patrick} Cattiaux }}
\address{{\bf {Patrick} CATTIAUX}\\ \'Ecole Polytechnique, CMAP, F- 91128 Palaiseau cedex
\\ and Universit\'e Paul Sabatier, Institut de Math\'ematiques, Laboratoire de
Statistiques et Probabilit\'es, 118 route de Narbonne, F-31062 Toulouse Cedex 09, France.}
\email{cattiaux@cmapx.polytechnique.fr}

\author[P. Collet]{\textbf{\quad {Pierre} Collet }}
\address{{\bf {Pierre} COLLET}\\ \'Ecole Polytechnique, CPHT, CNRS-UMR7644,  F- 91128
Palaiseau cedex, France.} \email{collet@cpht.polytechnique.fr}

\author[A. Lambert]{\textbf{\quad {Amaury} Lambert }}
\address{{\bf {Amaury} LAMBERT}\\  \, Laboratoire de Probabilit\'es et Mod\`eles Al\'eatoires,
UPMC Univ Paris 06,
4 Place Jussieu,
F-75252 Paris Cedex 05, France.}
\email{amaury.lambert@upmc.fr}

\author[S. Martinez]{\textbf{\quad\quad\quad\quad\quad {Servet} Mart\'inez }}
\address{{\bf {Servet} MART\'INEZ}\\  \, Universidad de Chile, Facultad de Ciencias F\'isicas y
Matem\'aticas, Departamento de Ingenier\'ia Matem\'atica y Centro de Modelamiento
Matem\'atico, Casilla 170-3, Correo 3, Santiago,
Chile.} \email{smartine@dim.uchile.cl}

 \author[S. M\'el\'eard]{\textbf{\quad {Sylvie} M\'el\'eard }}
\address{{\bf {Sylvie} M\'EL\'EARD}\\ \'Ecole Polytechnique, CMAP, F-
  91128 Palaiseau
cedex, France.} \email{meleard@cmapx.polytechnique.fr}

\author[J. San Martin]{\textbf{\quad {Jaime} San Mart\'in }}
\address{{\bf {Jaime} SAN MART\'IN}\\  \, Universidad de Chile, Facultad de Ciencias F\'isicas y
Matem\'aticas, Departamento de Ingenier\'ia Matem\'atica y Centro de Modelamiento
Matem\'atico, Casilla 170-3, Correo 3, Santiago,
Chile.} \email{jsanmart@dim.uchile.cl}

\begin{abstract}
In this paper, we study quasi-stationarity for a  large class of
Kolmogorov diffusions. The main novelty here is that we allow the drift to go
to $- \infty$ at the origin, and the diffusion to have an entrance
boundary at $+\infty$. These diffusions arise as images, by a
deterministic map, of generalized Feller diffusions, which
themselves are obtained as limits of rescaled birth--death
processes. Generalized Feller diffusions take nonnegative values
and are absorbed at zero in finite time with probability $1$.  An
important  example is the logistic Feller diffusion.

We give sufficient conditions on the drift near $0$ and near $+
\infty$ for the existence of quasi-stationary distributions, as
well as rate of convergence in the Yaglom limit and existence of the $Q$-process. We
also show that under these conditions, there is exactly one
quasi-stationary distribution, and that this distribution attracts all initial distributions
under the conditional evolution, if and only if $+\infty$ is an entrance
boundary. In particular this gives a sufficient condition for the uniqueness
of quasi-stationary distributions.
In the proofs spectral theory plays an important role on $L^2$ of the
reference measure for the killed
process.
\end{abstract}

\maketitle

\bigskip

\textit{ Key words.} quasi-stationary distribution, birth--death process, population dynamics, logistic growth, generalized Feller diffusion, Yaglom limit, convergence rate, $Q$-process, entrance boundary at infinity.
\bigskip

\textit{ MSC 2000 subject.} Primary 92D25; secondary 37A30, 60K35, 60J60, 60J85, 60J70.\\
\bigskip

\section{\bf Introduction}\label{Intro}

The main motivation of this work is the existence, uniqueness and domain of
attraction of  quasi-stationary distributions for some diffusion
models arising from population dynamics. After a change of
variable, the problem is stated in the framework of Kolmogorov
diffusion processes with a drift behaving like $-1/2x$ near the
origin. Hence, we shall study quasi-stationarity for the larger
class of one-dimensional Kolmogorov diffusions (drifted Brownian
motions), with drift possibly exploding at the origin.
\medskip

Consider a one-dimensional drifted Brownian motion on $(0,\infty)$
\begin{equation}\label{eqdiff}
dX_t \, = \, dB_t \, - \, q(X_t) \, dt \quad , \quad X_0=x > 0
\end{equation}
where $q$ is defined and  $C^1$  on $(0,\infty)$ and $(B_t; t\ge
0)$ is a standard one-dimensional Brownian motion. In particular
$q$ is allowed to explode at the origin. A pathwise unique
solution of \eqref{eqdiff} exists up to the explosion time
$\tau$.
We denote $T_y$  the first time the
process hits $y \in (0,\infty)$ (see \cite{IW} chapter VI section 3)
before the explosion
$$
T_y=\inf\{0\le t < \tau: X_t=y\}
$$
We denote by $T_\infty=\lim\limits_{n\to \infty} T_n$ and
$T_0=\lim\limits_{n\to \infty} T_{1/n}$.
Since $q$ is regular in $(0,\infty)$ then $\tau=T_0\wedge T_\infty$.

The law of the process starting from $X_0$ with distribution $\nu$
will be denoted by $\P_\nu$.
A \emph{quasi-stationary distribution} (in short {\it q.s.d.})
 for $X$ is a probability measure $\nu$
supported on $(0,\infty)$ satisfying for all $t\ge 0$
\begin{equation}
 \PP_{\nu}(X_t\in A \mid T_0>t) = \nu(A), \;\; \forall
 \mbox{ Borel set }A\subseteq (0,\infty).
\end{equation}
By definition a {\it q.s.d.} is a fixed point of the conditional
evolution. The \emph{Yaglom limit} $\pi$ is defined as the limit
in distribution
$$
\pi(\bullet)=\lim_{t\to \infty}\PP_x(X_t\in \bullet \mid T_0>t),
$$
provided this limit exists and is independent of the initial condition $x$.
The Yaglom limit is a {\it q.s.d.} (see Lemma \ref{qlid-qsd}).

We will also study the existence of the so-called $Q$-process which is obtained as the law of
the process $X$ conditioned to be never extinct, and it is defined as follows. For any $s\ge 0$ and
any Borel set ${B} \subseteq C([0,s])$ consider
$$
\Q_x (X\in {B}) = \lim_{t\to\infty} \P_x(X\in {B}\mid T_0>t).
$$
When it exists, this limit procedure defines the law of a diffusion that never reaches 0.
\medskip

The reason for studying such diffusion processes with a possibly exploding drift at the origin
comes from our interest in the following \textit{generalized Feller diffusion processes}
\begin{equation}
\label{eqn : generalized Feller i}
 dZ_t=\sqrt{\gamma Z_t}dB_t+h(Z_t)dt ,\; Z_0=z>0 ,
\end{equation}
where $h$ is a nice function satisfying $h(0)= 0$.

Notice that  $z=0$ is an absorbing state for  $Z$. This means that if $Z_0=0$ then
$Z_t=0$ for all $t$ is the unique solution of \eqref{eqn : generalized
Feller i} (see \cite{IW}).

If we define $X_t=2 \sqrt{Z_t/\gamma}$ then
\begin{equation}
\label{assogenediff}
dX_t = dB_t - \frac{1}{X_t} \, \left(\frac 12 -\frac{2}{\gamma} \,
h\left(\frac{\gamma X_t^2}{4}\right)\right) dt
\, \textrm{ , } \, X_0=x=2\sqrt{z/\gamma}>0 \, ,
\end{equation}
so that $X$ is a drifted Brownian motion as in \eqref{eqdiff} where  $q(x)$ behaves
like $1/2x$ near the origin. The process $Z$ is obtained after rescaling some sequences of
birth--death processes arising from population dynamics.

Of particular interest is the case $h(z)=rz - cz^2$ (logistic case), for which we obtain $$q(x) =
\frac{1}{2x} - \frac{r x}{2} + \frac{c \gamma x^3}{8} \, .$$

A complete description of these models is performed in the final section of the paper, where their
biological meaning is also discussed. Of course quasi-stationary distributions for $Z$ and $X$
are related by an immediate change of variables, so that the results on $X$ can be immediately
translated to results on $Z$.
\medskip

The study of quasi-stationarity is a long standing problem (see \cite{Pol} for a regularly
updated extensive bibliography and \cite{FKMP,Goss,SVJ} for the Markov chain case).
For Kolmogorov diffusions, the theory started with Mandl's paper
\cite{Mandl} in 1961, and was then developed by many authors (see in particular \cite{CMM,MSM,SE05}). All
these works  assume Mandl's conditions, which are not  satisfied in the situation
described above, since in particular the drift is not bounded near $0$. It is worth noticing that
the behavior of $q$ at infinity also may violate Mandl's conditions, since in the logistic case for
instance,
$$
\int_1^\infty e^{-Q(z)}\int_1^z e^{Q(y)} dy\, dz<\infty,
$$
where $Q(y):=2\int_1^y q(x)\, dx$.

This unusual situation prevents us from using earlier results on
{\it q.s.d.}'s of solutions of Kolmogorov equations. Hence we are
led to develop new techniques to cope with this
situation.

In Section \ref{sectiondiff} we start with the study of a general
Kolmogorov diffusion process on the half line and introduce the
hypothesis (H1) that ensures to reach 0 in finite time with
probability 1.  Then we introduce the measure $\mu$, not
necessarily finite, defined as 
$$
\mu(dy):=e^{-Q(y)}\, dy,
$$
 which is
the speed measure of $X$. We describe the Girsanov transform and
show how to use it in order to obtain $\L^2(\mu)$ estimates for
the heat kernel (Theorem \ref{thmloit}). The key is the following:
starting from any $x>0$, the law of the process at time $t$ is
absolutely continuous with respect to $\mu$ with a density
belonging to $\L^2(\mu)$ (and explicit bounds). In the present
paper we work in $\L^p(\mu)$ spaces rather than $\L^p(dx)$, since 
it greatly simplifies the presentation of the spectral theory.

This spectral theory is done in Section \ref{sectiongap}, where we introduce the hypothesis (H2):
$$
\lim_{x \to \infty} q^2(x) - q'(x)=\infty\,
\textrm{ , } \, C:=-\inf_{x \in (0,\infty)} q^2(x)-q'(x) <\infty.
$$
This hypothesis ensures the discreteness of the spectrum (Theorem \ref{thmspectrum}).
The ground state $\eta_1$ (eigenfunction associated to
the bottom of the spectrum) can be chosen nonnegative, even positive as we will see, 
and furnishes the natural candidate $\eta_1 d\mu$ for a {\it q.s.d.}. The only thing to check is
that $\eta_1 \in \L^1(\mu)$ which is not immediate since $\mu$ is possibly unbounded.

Section \ref{sec : eigenfunctions}  gives some sharper properties of the eigenfunctions defined in
the previous section, using in particular properties of the Dirichlet heat kernel.
More specifically, we introduce two independent hypotheses, (H3) and (H4), either of which ensures that
the eigenfunctions belong to $\L^1(\mu)$ (Propositions \ref{thmqsd} and \ref{propspt}). Hypothesis (H3) is
$$
\int_0^1 \frac{1}{q^2(y) - q'(y)+C+2}\,\mu(dy) < \infty,
$$
and (H4) is
$$
\int_1^\infty e^{-Q(y)}\, dy <\infty\, \textrm{ , } \,\int_0^1  ye^{-Q(y)/2}\, dy <\infty.
$$
Section \ref{sectionyaglom} contains the proofs of the
existence of the Yaglom limit (Theorem \ref{thm : Yaglom}) as well as the exponential decay
to equilibrium (Proposition \ref{ratecv}), under hypotheses
(H1) and (H2), together with either (H3) or (H4). Section \ref{secQ} contains the results on the
$Q$-process (Corollaries \ref{corq} and \ref{corjavaitort}).

In Section \ref{sectionattrac} we introduce condition (H5) which is equivalent to the
existence of an entrance law at $+\infty$, that is, the repelling force at infinity imposes to the
process starting from infinity to reach any finite interval in finite time.
The process is then said to ``come down from infinity''. We show
that the process comes down from infinity if and only if there
exists a unique {\it q.s.d.} which attracts any initial law under
the conditional evolution (Theorem \ref{thm : uniqueness}). In
particular this theorem gives sufficient conditions for uniqueness
of {\it q.s.d.}'s. In the context of birth and death chains the
equivalence between uniqueness of a {\it q.s.d.} and  ``come down
from infinity'' has been proved in \cite[Theorem 3.2]{VD}.


The final section contains the description of the underlying biological models,
as well as the application of the whole theory developed in the previous sections
to these models (Theorem \ref{thm : main result}).

In the following statement (which is basically Theorem \ref{thm : main result}), we record the main results of this paper in terms of the generalized Feller diffusion solution of \eqref{eqn : generalized Feller i}. We say that $h$ satisfies the condition (HH) if
$$
(i)\lim_{x \to \infty}\frac{h(x)}{\sqrt{x}}=-\infty,\qquad (ii)\lim_{x\to \infty}\frac{x
h'(x)}{h(x)^2} =0.
$$
\begin{theorem}
If $h$ satisfies (HH),
then for all initial laws with bounded support, the law of $Z_t$ conditioned on $\{Z_t\neq 0\}$
converges exponentially fast to a probability measure $\nu$, called the Yaglom limit.

The process $Z$ conditioned to be never extinct is well-defined and is called the $Q$-process. The $Q$-process converges in distribution, to
its unique invariant probability measure. This probability measure is absolutely continuous w.r.t.
$\nu$ with a nondecreasing Radon--Nikodym derivative.

If in addition,  the following integrability condition is satisfied
$$\int_1^{\infty}\frac{dx}{-h(x)}<\infty ,$$
then $Z$ comes down from infinity and the convergence of the conditional one-dimensional distributions
holds for all initial laws, so that the Yaglom limit $\nu$ is the
unique quasi-stationary distribution.
\end{theorem}

\section{\bf One dimensional diffusion processes on the positive
half line}
\label{sectiondiff}

Associated to $q$ we consider the functions
\begin{equation}\label{eqtest}
\Lambda(x)  =  \int_1^x \, e^{Q(y)} \, dy \qquad \textrm{and}\qquad
 \kappa(x)  =  \int_1^x \, e^{Q(y)} \,
\left(\int_1^y \, e^{- Q(z)} \, dz\right) \, dy \, ,
\end{equation}
where we recall that $Q(y) = \int_1^y \, 2 q(u) du$. Notice that
$\Lambda$ is the scale function for $X$.
\medskip

For most of the results in this paper we shall assume sure absorption at zero, that is
\begin{equation}\label{eqexplo}
\hspace{-5.5cm}\hbox{ {\bf Hypothesis (H1):}  for all } x>0\  ,\quad \P_x(\tau = T_0 <  T_\infty) = 1 \,
.
\end{equation}
It is well known (see e.g. \cite{IW} chapter VI Theorem 3.2) that
\eqref{eqexplo} holds if and only if
\begin{equation}\label{eqexplo2}
\Lambda(\infty)=\infty \, \textrm{ and } \, \kappa(0^+)< \infty \, .
\end{equation}

\medskip

We notice that (H1) can be written as
$\P_x(\lim\limits_{t \to \infty} X_{t\wedge \tau}=0)=1$.

\begin{example}\label{exdiff}
The main cases that we are interested in are the following ones.

\begin{enumerate}

\item[(1)] When $X$ is defined by \eqref{assogenediff} associated to
the generalized Feller diffusion $Z$.
It is direct to show that $Q(x)$ behaves like $\log(x)$ near $0$ hence $\kappa(0^+)<\infty$.
The logistic Feller diffusion corresponds to $h(z)=rz -c z^2$ for some constants $c$ and $r$.
It is easily seen that \eqref{eqexplo2} is satisfied in this case
provided $c>0$ or $c=0$ and $r<0$.

\item[(2)] When the drift is bounded near $0$, in which case $\kappa(0^+)<\infty$.

\hfill{$\diamondsuit$}
\end{enumerate}
\end{example}

We shall now discuss some properties of the law of $X$ up to
$T_0$. The first result is a Girsanov type result.

\begin{proposition}
\label{propgirsanov} Assume (H1). For any Borel bounded function
$F:C([0,t],(0,\infty))\to \RR$ it holds
$$\E_x \left[F(X) \, \BBone_{t<T_0}\right] \, =
\, \E^{\W_x} \left[ F(\omega) \, \BBone_{t<T_0(\omega)} \, \exp
\left(\frac 12 \, Q(x) -  \frac 12 \, Q(\omega_t) - \, \frac 12 \,
\int_0^t \, (q^2 - q')(\omega_s) ds\right)\right]$$ where
$\E^{\W_x}$ denotes the expectation w.r.t. the Wiener measure
starting from $x$, and $\E_x$ denotes the expectation with respect
to the law of $X$ starting also from $x$.
\end{proposition}
\begin{proof}
It is enough to show the result for $F$ nonnegative and bounded.
Let $x>0$ and consider  $\varepsilon \in (0,1)$ such that
$\varepsilon \leq x \leq 1/\varepsilon$. Also we define
$\tau_\varepsilon=T_\varepsilon \wedge T_{1/\varepsilon}$. Choose
some $\psi_\varepsilon$ which is a nonnegative $C^\infty$ function
with compact support, included in $]\varepsilon/2, 2/\varepsilon[$
such that $\psi_\varepsilon(u)=1$ if $\varepsilon \leq u \leq
1/\varepsilon$. The law of the diffusion \eqref{eqdiff} coincides
up to $\tau_\varepsilon$ with the law of a similar diffusion
process $X^\varepsilon$ obtained by replacing $q$ with the cutoff
$ q_{\varepsilon}=q\psi_\varepsilon$. For the latter we may apply
the Novikov criterion ensuring that the law of $X^\varepsilon$ is
given via the Girsanov formula. Hence
\begin{eqnarray*}
\E_x \left[F(X) \, \BBone_{t<\tau_\varepsilon}\right] & = &
\E^{W_x} \left[ F(\omega) \, \BBone_{t<\tau_\varepsilon(\omega)}
\, \exp \left(\int_0^t \, - q_\varepsilon(\omega_s) d\omega_s - \,
\frac
12 \, \int_0^t \, (q_\varepsilon)^2(\omega_s) ds\right)\right] \\
& = &  \E^{W_x} \left[ F(\omega) \,
\BBone_{t<\tau_\varepsilon(\omega)} \, \exp \left(\int_0^t \, -
q(\omega_s) d\omega_s - \, \frac 12 \, \int_0^t \, q^2(\omega_s)
ds\right)\right]\\& = & \E^{W_x} \left[ F(\omega) \,
\BBone_{t<\tau_\varepsilon(\omega)} \, \exp \left(\frac 12 \, Q(x)
- \frac 12 \, Q(\omega_t) - \, \frac 12 \, \int_0^t \, (q^2 -
q')(\omega_s) ds\right)\right]\;.
\end{eqnarray*}
The last equality is obtained integrating by parts the stochastic integral. But
$\BBone_{t<\tau_\varepsilon}$ is non-decreasing in $\varepsilon$
and converges  almost surely to $\BBone_{t<T_0}$ both for $\P_x$
(thanks to (H1)) and $\W_x$.
It remains to use Lebesgue monotone convergence theorem
to finish the proof.
\end{proof}

The next theorem is inspired by the calculation in Theorem 3.2.7
of \cite{Ro99}. It will be useful to introduce the following
measure defined on $(0,\infty)$
\begin{equation}
\label{eqmu}
\mu(dy) := e^{-Q(y)} \, dy \, .
\end{equation}
Note that $\mu$ is not necessarily finite.

\begin{theorem}
\label{thmloit} Assume (H1). For all $x>0$ and all $t>0$ there
exists some density $r(t,x,.)$ that satisfies
$$\E_x[f(X_t) \, \BBone_{t<T_0}] = \int_0^{\infty} \, f(y) \,
r(t,x,y) \, \mu(dy)$$
for all bounded Borel $f$.

If in addition there exists some $C>0$ such that $q^2(y)-q'(y)
\geq - C$ for all $y>0$, then for all $t>0$ and all $x>0$,
$$\int_0^{\infty} \, r^2(t,x,y) \, \mu(dy) \, \leq (1/2\pi t)^{\frac 12}
\, e^{Ct} \, e^{Q(x)} \, .$$
\end{theorem}
\begin{proof}
Define
$$G(\omega) = \BBone_{t<T_0(\omega)} \,
\exp \left(\frac 12 \, Q(\omega_0) - \frac 12 \, Q(\omega_t) - \,
\frac 12 \, \int_0^t \, (q^2 - q')(\omega_s) ds\right) \, .$$
Denote by
$$e^{-v(t,x,y)}= (2\pi t)^{- \frac 12} \, \exp \left( -
\frac{(x-y)^2}{2t}\right)$$ the density at time $t$ of the
Brownian motion starting from $x$. According to Proposition
\ref{propgirsanov} we have
\begin{eqnarray*}
\E_x[f(X_t) \, \BBone_{t<T_0}] & = & \E^{\W_x}[f(\omega_t)
\E^{\W_x}[G|\omega_t]] \\ & = & \int \, f(y) \,
\E^{\W_x}[G|\omega_t=y] \, e^{-v(t,x,y)} \, dy \\ & = &
\int_0^{\infty} \, f(y) \, \E^{\W_x}[G|\omega_t=y] \,
e^{-v(t,x,y)+Q(y)} \, \mu(dy) \, ,
\end{eqnarray*}
because $\E^{\W_x}[G|\omega_t=y]=0$ if $y\leq 0$. In other words,
the law of $X_t$
restricted to non-extinction has a density with respect to $\mu$
given by
$$
r(t,x,y)=\E^{\W_x}[G|\omega_t=y] \, e^{-v(t,x,y)+Q(y)} \, .
$$
Hence,
\begin{eqnarray*}
\int_0^{\infty} \, r^2(t,x,y) \, \mu(dy) & = & \int \,
\left(\E^{\W_x}[G|\omega_t=y] \, e^{-v(t,x,y)+Q(y)}\right)^2 \,
e^{-Q(y)+v(t,x,y)} \, e^{-v(t,x,y)} \, dy \\ & = & \E^{\W_x}
\left[e^{-v(t,x,\omega_t)+Q(\omega_t)} \,
\left(\E^{\W_x}[G|\omega_t]\right)^2 \right] \\ & \leq & \E^{\W_x}
\left[e^{-v(t,x,\omega_t)+Q(\omega_t)} \, \E^{\W_x}[G^2|\omega_t]
\right] \\ & \leq & e^{Q(x)} \, \E^{\W_x}
\left[\BBone_{t<T_0(\omega)} \, e^{-v(t,x,\omega_t)} \, e^{-
\int_0^t \, (q^2 - q')(\omega_s) ds}\right] \, ,
\end{eqnarray*}
where we have used Cauchy-Schwarz's inequality. Since $e^{-v(t,x,.)} \leq
(1/2\pi t)^{\frac 12}$ the proof is completed.
\end{proof}

\begin{remark}
\label{remcondition} It is interesting to discuss a little bit the
conditions we have introduced.

\begin{enumerate}
\item[(1)] Since $q$ is assumed to be regular,
the condition $q^2 - q'$ bounded from below has to be checked only near
infinity or near $0$.

\item[(2)] Consider the behavior near infinity. Let us show
that if $\liminf_{y \to \infty} (q^2(y)-q'(y)) =
-\infty$ then $\limsup_{y \to \infty} (q^2(y)-q'(y)) > -\infty$
i.e. the drift $q$ is strongly oscillating. Indeed, assume that $q^2(y)-q'(y)
\to -\infty$ as $y \to \infty$. It follows that $q'(y) \to
\infty$, hence $q(y) \to \infty$. For $y$ large enough we may
thus write $q(y)=e^{u(y)}$ for some $u$ going to infinity at
infinity. So $e^{2u(y)} (1-u'(y)e^{-u(y)}) \to - \infty$ implying
that $u' e^{-u} \geq 1$ near infinity. Thus if $g = e^{-u}$ we
have $g'\leq -1$ i.e. $g(y) \to - \infty$ as $y \to \infty$ which
is impossible since $g$ is nonnegative.

\item[(3)] If $X$ is given by (\ref{assogenediff}) we have
$$
q(y)=\frac 1y \left(\frac 12 \, - \, \frac{2}{\gamma}
h\left(\frac{\gamma y^2}{4}\right)\right)\,.
$$
Hence, since $h$ is of class $C^1$ and $h(0)=0$, $q^2(y)-q'(y)$ behaves near $0$
like $\frac{3}{4y^2} $ so that $q^2 - q'$ is bounded
from below near 0 (see Appendix for further conditions fulfilled by $h$ to get
the same result near $\infty$).
\hfill{$\diamondsuit$}
\end{enumerate}
\end{remark}

\bigskip

\section{\bf $\L^2$ and spectral theory of the diffusion process}
\label{sectiongap}

Theorem \ref{thmloit} shows that for a large family of initial laws, the
distribution of $X_t$ before extinction has a density belonging to
$\L^2(\mu)$. The measure $\mu$ is natural since the kernel of the
killed process is symmetric in $\L^2(\mu)$, which allows us to use spectral theory.

Let $C^\infty_0((0,\infty))$ be the vector space of infinitely differentiable functions on $(0,\infty)$
with compact support.
We denote
$$\langle f,g \rangle_\mu \, = \, \int_0^{\infty} \,  f(u) g(u) \mu(du) \, .$$
Consider the symmetric form
\begin{equation}\label{eqform}
\mathcal E (f,g)=\langle f', g' \rangle_\mu \, \quad\mathrm{ , }\quad\,D(\mathcal E)=C^\infty_0((0,\infty)).
\end{equation}
This form is Markovian and closable. The proof of the latter
assertion is similar to the one of Theorem 2.1.4 in \cite{Fuku}
just replacing the real line by the positive half line. Its
smallest closed extension, again denoted by $\mathcal E$, is thus
a Dirichlet form which is actually regular and local. According to
the theory of Dirichlet forms (see \cite{Fuku} or \cite{FOT}) we
thus know that
\begin{itemize}
\item there exists a non-positive self adjoint operator $L$
on $\L^2(\mu)$ with domain $D(L)\supseteq C^\infty_0((0,\infty) $
such that for all $f$ and $g$ in
$C^\infty_0((0,\infty))$ the following holds (see \cite{Fuku}
Theorem 1.3.1)
\begin{equation}
\label{eqdirichlet}
\mathcal E (f,g) = - \, 2 \, \int_0^{\infty} \, f(u) \, Lg(u) \,
\mu(du) \, = \, - \, 2 \, \langle f,Lg \rangle_\mu \, .
\end{equation}
We point out that for $g\in
C^\infty_0((0,\infty))$,
$$
Lg=\frac 12 \, g'' - q g'.
$$
\item $L$ is the generator
of a strongly continuous symmetric semigroup of contractions on
$\L^2(\mu)$ denoted by $(P_t)_{t\geq 0}$. This semigroup is
(sub)-Markovian, i.e. $0\leq P_t f \leq 1$ $\mu$ a.e. if $0 \leq f
\leq 1$ (see \cite{Fuku} Theorem 1.4.1).
\item There exists
a unique $\mu$-symmetric Hunt process with continuous sample
paths (i.e. a diffusion process) up to its explosion time $\tau$
whose Dirichlet form is $\mathcal E$ (see \cite{Fuku} Theorem
6.2.2)
\end{itemize}
The last assertion implies that, for $\mu$ quasi all
$x>0$ (that is, except for a set of zero capacity, see \cite{Fuku} for details),
one can find a probability measure $\Q_x$ on
$C(\R^+,(0,\infty))$ such that for all $f\in
C^\infty_0((0,\infty))$,
$$f(\omega_{t\wedge \tau})-f(x)-\int_0^{t \wedge \tau} \, Lf(\omega_s) ds$$
is a local martingale with quadratic variation  $\int_0^{t \wedge \tau}
|f'|^2(\omega_s) ds$. Due to our hypothesis $q\in C^{1}(0,\infty)$
we know that this martingale problem admits a unique solution 
(see for example \cite{Jacod} page 444). On the other hand,
using It\^o's formula we know that under $\P_x$, 
the law of $(X_{t\wedge \tau})$ is also a solution to this martingale
problem. 

\medskip

The conclusion is that  the semigroup $P_t $ and the semigroup 
induced by the strong Markov process $(X_{t\wedge \tau})$ 
coincide on the set of smooth and compactly supported functions.
Therefore, for all  $f\in \L^2(\mu)$ we have that
$$
P_t f(x)=\EE[f(X_t) \BBone_{t<\tau}].
$$

\medskip

Let $(E_\lambda: \lambda \ge 0)$ be the spectral family of $-L$.  We can
restrict ourselves to the case $\lambda\ge 0$ because $-L$ is nonnegative. Then $\forall
\, t\geq 0, \, f, \, g\in \L^2(\mu)$,
\begin{equation}
\label{eqspectral} \int P_t f \, g \, d\mu = \int_0^{\infty} \,
e^{-\lambda t} \, d\langle E_\lambda f,g\rangle_\mu \, .
\end{equation}

We notice that if absorption is sure, that is (H1) holds, this semigroup coincides with
the semigroup of $X$ killed at $0$, that is
$P_t f(x)=\E[f(X_t^x) \, \BBone_{t<T_0}]$.

Note that for $f\in \L^2(\mu)$ and all closed interval
$K \subset (0,\infty)$,
\begin{eqnarray*}
\int (P_t f)^2 d\mu & = & \int (P_t (f \BBone_K + f \BBone_{K^c}))^2 d\mu
\\ & \leq & 2 \int (P_t (f \BBone_K))^2 d\mu +
2 \int (P_t (f \BBone_{K^c}))^2 d\mu \\
& \leq & 2 \int (P_t (f \BBone_K))^2 d\mu + 2 \int (f
\BBone_{K^c})^2 d\mu \, .
\end{eqnarray*}
We may choose $K$ large enough in order that the second term in
the latter sum is bounded by $\varepsilon$. Similarly we may
approximate $f \BBone_K$ in $\L^2(\mu)$ by $\tilde f \BBone_K$ for
some continuous and bounded $\tilde f$, up to $\varepsilon$
(uniformly in $t$). Now, thanks to (H1) we know that $P_t
(\tilde f \BBone_K)(x)$ goes to $0$ as $t$ goes to infinity for
any $x$. Since
$$\int (P_t (\tilde f \BBone_K))^2 d\mu =
\int_K \tilde f \, P_{2t} (\tilde f \BBone_K) \, d\mu \,,$$ we may
apply Lebesgue bounded convergence theorem and conclude that $\int
(P_t (\tilde f \BBone_K))^2 d\mu \to 0$ as $t \to \infty$. Hence,
we have shown that,
\begin{equation}
\label{eqnozero} \forall \, f \in \L^2(\mu) \quad \int
(P_t f)^2 d\mu \, \to \, 0 \, \textrm{ as } \, t \, \to \, \infty
\, .
\end{equation}

Now we shall introduce the main assumption on $q$ for the spectral aspect
of the study.

\begin{equation}
\label{defspectrum}
\hspace{-0.05cm}\hbox{{\bf Hypothesis (H2)} }\; C=-\inf\limits_{y\in (0,\infty)} q^2(y)-q'(y)<\infty
\hbox{ and } \lim_{y \to \infty} \, q^2(y) - q'(y) = + \infty \, .
\end{equation}

\begin{proposition}
Under (H2), $|q(x)|$ tends to infinity as $x \to \infty$, and $q^-(x)$ or
$q^+(x)$ tend to 0 as $x\downarrow 0$. If in addition (H1)
holds then $q(x)\to \infty$, as $x\to \infty$.
\end{proposition}
\begin{proof}

Since $q^2-q'$ tends to $\infty$  as $x\to \infty$, $q$ does not
change sign for large $x$. If $q$ is bounded near
infinity we arrive to a contradiction because $q'$ tends to
$-\infty$ and therefore $q$ tends to $-\infty$ as well. So $q$ is unbounded.
If $\liminf\limits_{x\to \infty} |q(x)|=a<\infty$ then we can
construct a sequence $x_n\to \infty$ of local maxima, or local
minima of $q$ whose value $|q(x_n)|<a+1$, but then
$q^2(x_n)-q'(x_n)$ stays bounded, which is a contradiction.

Now we prove that $q^-(x)$ or $q^+(x)$ tend to 0 as $x\downarrow
0$. In fact, assume there exist an $\epsilon>0$ and a sequence $(x_n)$
with $0<x_n \downarrow 0$ such that $q(x_{2n})=-\epsilon,
q(x_{2n+1})=\epsilon$. Then we can construct another sequence
$z_n\downarrow 0$ such that $|q(z_n)|\le \epsilon$ and $q'(z_n)\to
\infty$, contradicting (H2).

Finally, assume (H1) holds.  If $q(x)\le -1$ for all $x>x_0$ we arrive to a contradiction.
Indeed, for all $t$
$$
\P_{x_0+1} (T_0>t)\ge \P_{x_0+1} (T_{x_0} > t)\ge \P_{x_0+1} (T_{x_0} =
\infty) \;.
$$
The assumption  $q(x)\le -1$ implies that $X_{t}\ge B_{t}+t$ while $t\le T_{x_0}$, and therefore
$$
\P_{x_0+1} (T_0>t)\ge \P_{x_0+1}\left(B_{t}+t\;\text{hits}\;\infty\;\text{before}\;x_{0}
\right)=\frac{e^{-2(x_{0}+1)}-e^{2x_{0}}}{e^{-2\infty}-e^{2x_{0}}}
=1-e^{-2}
$$
where we have used that $\big(\exp(-2(B_{t}+t))\big)$ is a martingale.
This contradicts (H1) and we have $q(x)\to \infty$ as $x\to \infty$.
\end{proof}
\medskip

We may now state the following result.

\begin{theorem}\label{thmspectrum}
If (H2) is satisfied, $-L$ has a purely discrete spectrum $0\le
\lambda_1 < \lambda_2 < ...$. Furthermore each $\lambda_i$
($i\in\mathbb{N}$)
 is
associated to a unique (up to a multiplicative constant)
eigenfunction $\eta_i$ of class $C^2((0,\infty))$, which also satisfies
the ODE
\begin{equation}
\label{ODE}
\frac12\eta_i''-q \eta_i'=-\lambda_i \eta_i.
\end{equation}
The sequence $(\eta_i)_{i\geq 1}$ is an orthonormal basis of
$~\L^2(\mu)$, $\eta_1$ can be chosen to be strictly positive in
$(0,\infty)$.

For $g \in \L^2(\mu)$,
$$P_tg=\sum_{i\in\mathbb{N}} e^{-\lambda_it}\langle \eta_i,g\rangle_\mu \eta_i
\quad {\rm in} \quad \L^2(\mu),$$  then for  $f,g \in
\L^2(\mu)$,
$$
\lim_{t \to \infty} e^{\lambda_1 t} \langle  g,P_t f \rangle_\mu
=  \langle \eta_1, f\rangle_\mu \, \langle \eta_1, g\rangle_\mu \,
.
$$
 If, in addition, (H1) holds, then
$\lambda_1> 0$.
\end{theorem}

\begin{proof}
For $f \in \L^2(dx)$,  define $\tilde{P_t}(f)=e^{-Q/2} P_t(f \,
e^{Q/2})$,  which exists in  $\L^2(dx)$ since $f \, e^{Q/2} \in
\L^2(\mu)$. $(\tilde{P_t})_{t\geq 0}$ is then a strongly
continuous semigroup in $\L^2(dx)$, whose generator $\tilde L$
coincides on $C^\infty_0((0,\infty))$ with $\frac 12
\frac{d^2}{dx^2} - \frac 12 \, (q^2 - q')$ since
$C^{2}_{0}(0,\infty)\subset D(L)$, and $e^{Q/2}\in
C^{2}(0,\infty)$. The spectral theory of such a Schr\"{o}dinger
operator on the line (or the half line) is well known, but here
the potential $v=(q^2 - q')/2$ does not necessarily belong to
$\L_{loc}^\infty$ near $0$ as it is generally assumed. We shall
use \cite{BS} chapter 2.

First we follow the proof of Theorem 3.1 in \cite{BS}. Since we
have assumed that $v$ is bounded from below by $-C/2$, we may
consider $H=\tilde L - (C/2+1)$, i.e. replace $v$ by $v+C/2+1=w\geq
1$, hence translate the spectrum. Since for $f\in C_{0}^{\infty}(0,\infty)$
\begin{equation}\label{eqschro}
-  (Hf,f) := - \int_0^{\infty} \!\! Hf(u) f(u) \, du =
\int_0^{\infty}\!\! \left(|f'(u)|^2/2+ w(u) f^2(u)\right) du
\geq \int_0^{\infty} f^2(u) \, du,
\end{equation}
 $H$ has a bounded inverse operator. Hence the spectrum of $H$
(and then the one  of $\tilde L$) will be discrete as soon as $H^{-1}$ is a
compact operator, i.e. as soon as $M=\{f\in D(H) \, ; \, - (Hf,f)
\leq 1\}$ is relatively compact. This is shown in \cite{BS} when
$w$ is locally bounded, in particular bounded near 0. If $w$ goes
to infinity at 0, the situation is even better since our set $M$
is included into the corresponding one with $w\approx 1$ near the
origin, which is relatively compact thanks to the asymptotic
behavior of $v$. The conclusion of Theorem 3.1 in \cite{BS} is
thus still true in our situation, i.e. the spectrum is discrete.

The discussion in Section 2.3 of \cite{BS}, pp. 59-69, is only
concerned with the asymptotic behavior (near infinity) of the
solutions of $f'' - 2wf=0$. Nevertheless, the results there applies to
our case.  All eigenvalues of $\tilde L$ are thus simple (Proposition 3.3 in
\cite{BS}), and of course the corresponding set of normalized eigenfunctions
$(\psi_k)_{k \geq 1}$ is an orthonormal basis of $\L^2(dx)$.

The system $(e^{Q/2} \, \psi_k)_{k \geq 1}$ is thus an orthonormal
basis of $\L^2(\mu)$, each $\eta_k=e^{Q/2} \, \psi_k$ being an
eigenfunction of $L$. We can choose them to be $C^2((0,\infty))$ and
they satisfy \eqref{ODE}.

For every $t>0$, and for every $g,f \in \L^2(\mu)$ we have
$$
\sum\limits_{k=1}^\infty e^{-\lambda_k t} \langle \eta_k, g
\rangle_\mu \langle \eta_k, f \rangle_\mu =\langle g, P_t f\rangle_\mu\;.
$$

In addition if $g$ and $f$ are nonnegative we get
$$
0 \leq \lim_{t \to \infty} e^{\lambda_1 t} \langle  g,P_t f \rangle_\mu
=  \langle \eta_1, f\rangle_\mu \, \langle \eta_1, g\rangle_\mu \,,
$$
since $\lambda_{1}<\lambda_{2}\le\ldots$ and the sum $\sum_{k=1}^\infty |\langle \eta_k, g \rangle_\mu \langle \eta_k,
f \rangle_\mu|$ is finite. It follows that $\langle \eta_1,
f\rangle_\mu $ and $\langle \eta_1, g\rangle_\mu$ have the same
sign. Changing $\eta_1$ into $- \eta_1$ if necessary, we may
assume that $\langle \eta_1, f\rangle_\mu\geq 0$ for any
nonnegative $f$, hence $\eta_1 \geq 0$. Since $P_t \eta_1
(x)=e^{-\lambda_1 t} \eta_1(x)$ and $\eta_1$ is continuous and not
trivial, we deduce that $\eta_1(x)>0$ for all $x>0$.

Since $L$ is non-positive, $\lambda_1\geq 0$. Now assume that (H1)
 holds. Using \eqref{eqnozero} we get for $g \in \L^2(\mu)$
$$
0=\lim_{t \to \infty} \langle P_t g,P_t g\rangle_\mu = \lim_{t
\to \infty} e^{- 2 \lambda_1 t} \langle g,\eta_1\rangle_\mu^2 \;,
$$
 showing that $\lambda_1 > 0$.
\end{proof}
\medskip

Moreover, we are   able to obtain a pointwise representation of the
density $r$.

\medskip

\begin{proposition}
\label{propdensity} Under (H1) and (H2) we have
\begin{equation}
\label{serie} r(t,x,y)=\sum\limits_{k=1}^\infty e^{-\lambda_k t}
\eta_k(x) \eta_k(y),
\end{equation}
uniformly on compact sets of
$(0,\infty)\times(0,\infty)\times(0,\infty)$.
\medskip

Therefore on compact sets of $(0,\infty)\times(0,\infty)$ we get
\begin{equation}
\label{firstlimit} \lim\limits_{t\to \infty} e^{\lambda_1 t}
r(t,x,y)=\eta_1(x) \eta_1(y).
\end{equation}

\end{proposition}
\begin{proof}
Using Theorems \ref{thmloit} and \ref{thmspectrum}, for every smooth
function $g$ compactly supported on $(0,\infty)$ we have
$$
\sum\limits_{k=1}^n e^{-\lambda_k t} \langle \eta_k, g \rangle_\mu^2
\le \sum\limits_{k=1}^\infty e^{-\lambda_k t} \langle \eta_k, g
\rangle_\mu^2 = \int\!\!\!\int g(x) g(y) r(t,x,y) e^{-Q(x) - Q(y)} dx dy.
$$
Then using the regularity of $\eta_k$ and $r$ we obtain, by letting
$g(y)dy$ tend to the Dirac measure in $x$, that
$$
\sum\limits_{k=1}^n e^{-\lambda_k t} \eta_k(x)^2  \le  r(t,x,x).
$$
Thus, the series $\sum\limits_{k=1}^\infty
e^{-\lambda_k t} \eta_k(x)^2$ converges pointwise, which by
Cauchy-Schwarz inequality implies the pointwise absolute
convergence of $\zeta(t,x,y):=\sum\limits_{k=1}^\infty
e^{-\lambda_k t} \eta_k(x) \eta_k(y)$ and the bound for all $n$
$$
\sum\limits_{k=1}^n e^{-\lambda_k t} |\eta_k(x) \eta_k(y)|\le
\sqrt{r(t,x,x)} \sqrt{r(t,y,y)}.
$$
Using Harnack inequality (see for example \cite{kry}) we get
$$
\sqrt{r(t,x,x)} \sqrt{r(t,y,y)}\le C_{K }r(t,x,y)
$$
for any $x$ and $y$ in the compact subset $K$ of $(0,\infty)$.
Using the dominated convergence theorem we obtain that for all
Borel functions $g,f$ with compact support in $(0,\infty)$
$$
\int\int g(x) f(y) \zeta(t,x,y) e^{-Q(x)-Q(y)} dx dy= \int\int
g(x) f(y) r(t,x,y) e^{-Q(x)-Q(y)} dx dy.
$$
Therefore $\zeta(t,x,y)=r(t,x,y)\; dxdy$-a.s., which proves the almost sure
version of (\ref{serie}).

\smallskip

Since $\eta_k$ are smooth eigenfunctions we get the pointwise
equality
$$
\begin{array}{ll}
e^{-\lambda_k t} \eta_k(x)^2 &= 
e^{-\lambda_k t/3} \langle r(t/3,x,\bullet), \eta_k\rangle_\mu 
\langle r(t/3,x,\bullet), \eta_k\rangle_\mu\\
&=\int\int r(t/3,x,y) r(t/3,x,z)
e^{-\lambda_k t/3} \eta_k(y) \eta_k(z) e^{-Q(z)-Q(y)} dy dz,
\end{array}
$$
which together with the fact $r(t/3,x,\bullet)\in \L^2(\mu)$ and
Theorem \ref{thmloit} allow us to
deduce
$$
\begin{array}{ll}
\sum_{k=1}^\infty e^{-\lambda_k t} \eta_k(x)^2 &=\int\int r(t/3,x,y) r(t/3,x,z)
\sum_{k=1}^\infty e^{-\lambda_k t/3} \eta_k(y) \eta_k(z) e^{-Q(z)-Q(y)} dy dz\\
&=\int\int r(t/3,x,y) r(t/3,x,z) r(t/3,y,z) e^{-Q(z)-Q(y)} dy dz=r(t,x,x).
\end{array}
$$
Dini's theorem then proves the uniform convergence in compacts of $(0,\infty)$
for the series
$$
\sum_{k=1}^\infty e^{-\lambda_k t} \eta_k(x)^2=r(t,x,x).
$$
By the Cauchy-Schwarz inequality we have for any $n$
$$
\left|\sum\limits_{k=n}^\infty e^{-\lambda_k t}
\eta_k(x) \eta_k(y)\right|
\le \left(\sum_{k=n}^\infty e^{-\lambda_k t}\eta_k(x)^2
\right)^{1/2}
 \left(\sum_{k=n}^\infty e^{-\lambda_k t}\eta_k(y)^2
\right)^{1/2} .
$$

This together with the dominated convergence
theorem yields \eqref{firstlimit}.
\end{proof}

In the previous theorem, notice that $\sum_k e^{-\lambda_k t}=\int
r(t,x,x) e^{-Q(x)} dx$ is the $\L^1(\mu)$ norm of $x\mapsto
r(t,x,x)$. This is finite if and only if $P_t$ is a trace-class operator on
$\L^2(\mu)$.

\section{\bf Properties of the eigenfunctions}
\label{sec : eigenfunctions}

In this section, we  study some properties of the eigenfunctions
$\eta_i$, including their  integrability with respect to $\mu$.

\begin{proposition}
\label{prop : whyme}
Assume that (H1) and (H2) are satisfied. Then $\int_1^\infty
\eta_1 e^{-Q} dx <\infty$, $F(x)=\eta_1'(x)e^{-Q(x)}$ is a
nonnegative decreasing function and the following limits exist
$$
F(0^+)=\lim\limits_{x\downarrow 0} \eta_1'(x)e^{-Q(x)} \in (0,\infty],\;\;
F(\infty)=\lim\limits_{x\to \infty} \eta_1'(x)e^{-Q(x)} \in [0,\infty).
$$
Moreover $\int\limits_0^\infty \eta_1(x) e^{-Q(x)} dx=\frac{F(0^+)-F(\infty)}{2\lambda_1}$.
In particular
$$
\eta_1 \in \L^1(\mu) \hbox{ if and only if } F(0^+)<\infty.
$$
The function $\eta_{1}$ is increasing and $\int\limits_1^\infty
e^{-Q(y)} dy<\infty$.
\end{proposition}

\begin{remark}
Note that $g=\eta_1 e^{-Q}$ satisfies the adjoint equation
$\frac12 g''+(qg)'=-\lambda_1 g$, and then $F(x)=g'(x)+2q(x)g(x)$ represents the flux
at $x$. Then $\eta_1\in \L^1(\mu)$ or equivalently $g \in \L^1(dx)$ if and only if the
flux at 0 is finite.
\hfill{$\diamondsuit$}
\end{remark}

\begin{proof}
Since $\eta_1$ satisfies $\eta_1''(x)-2q\eta_1'(x)=-2\lambda_1
\eta_1(x)$,  we obtain for $x_{0}$ and $x$ in $(0,\infty)$
\begin{equation}
\label{ecuacion}
\eta_1'(x)e^{-Q(x)}=\eta_1'(x_0)e^{-Q(x_0)}
-2\lambda_1 \int_{x_0}^x \eta_1(y) e^{-Q(y)} dy,
\end{equation}
and $F=\eta_1'e^{-Q}$ is decreasing.
Integrating further gives
$$
\eta_1(x)=\eta_1(x_0)+\int_{x_0}^x \left(\eta_1'(x_0)e^{-Q(x_0)}
-2\lambda_1 \int_{x_0}^z \eta_1(y) e^{-Q(y)} dy\right)
e^{Q(z)} dz.
$$
If for some $z_0>x_0$ it holds that
$\eta_1'(x_0)e^{-Q(x_0)}-2\lambda_1\int_{x_0}^{z_0} \eta_1(y) e^{-Q(y)} dy<0$,
then this inequality holds for all $z>z_0$ since the quantity
$$
\eta_1'(x_0)e^{-Q(x_0)}
-2\lambda_1 \int_{x_0}^z \eta_1(y) e^{-Q(y)} dy
$$
is decreasing in $z$.
 This implies that for large $x$ the function $\eta_1$ is negative, because $e^{Q(z)}$ tends to
$\infty$ as $z\to \infty$. This is a contradiction and we deduce that for all $x>0$
$$
2\lambda_1\int_x^\infty \eta_1(y) e^{-Q(y)} dy\le \eta_1'(x) e^{-Q(x)}.
$$
This implies that $\eta_1$ is increasing and, being nonnegative,
it is bounded near 0. In particular, $\eta_1(0^+)$ exists. Also we
deduce that $F\ge 0$ and that $\int\limits_1^\infty e^{-Q(y)}
dy<\infty$. We can take the limit as $x\to \infty$ in
\eqref{ecuacion} to get
$$
F(\infty)=\lim\limits_{x \to \infty} \eta_1'(x) e^{-Q(x)} \in [0,\infty),
$$
and $\eta'_1(x_0) e^{-Q(x_0)}=F(\infty)+2\lambda_1
\int_{x_0}^\infty \eta_1(y) e^{-Q(y)} dy$. From this equality the result
follows.
\end{proof}

In the next results we give some
sufficient conditions, in terms of $q$,
for the integrability of the eigenfunctions.
A first useful condition is the following one
$$
\hspace{-5.2cm} \textbf{Hypothesis (H3): }\,  \int_0^1 \, \frac{1}{q^2(y) - q'(y) +
C + 2} \, e^{-Q(y)}dy < \infty,
$$
where as before $C=-\inf\limits_{x>0} (q^2(x) - q'(x))$.

\begin{proposition}\label{thmqsd}
Assume that (H1), (H2) and (H3) are satisfied. Then $\eta_i$ belongs to
$\L^1(\mu)$ for all $i$.
\end{proposition}
\begin{proof}
Recall that $\psi_i=e^{-Q/2} \eta_i$ is an eigenfunction of the
Schr\"{o}dinger operator $H$ introduced in the proof of Theorem
\ref{thmspectrum}. Replacing $f$ by $\psi_i$ in \eqref{eqschro}
thus yields
$$
(C/2+1+\lambda_i) \int_0^{\infty} \psi_i^2(y)dy=\int_0^{\infty} (|\psi_i'|^2(y)/2+w(y)
\psi_i^2(y)) dy\;.
$$
Since the left hand side is finite, the right hand side is finite,
in particular
$$
\int_0^{\infty} w(y) \eta_i^2(y) \mu(dy) \, = \,
\int_0^{\infty} w(y) \psi_i^2(y) dy < \infty \, .
$$
 As a
consequence, using Cauchy-Schwarz inequality we get on one hand
\begin{eqnarray*}
\int_0^1 |\eta_i(y)| \mu(dy) & \leq & \left(\int_0^1 w(y) \,
\eta_i^2(y) \mu(dy)\right)^{\frac 12} \, \left(\int_0^1
\frac{1}{w(y)} \,  \mu(dy)\right)^{\frac 12} < \infty
\end{eqnarray*}
thanks to (H3). On the other hand
\begin{eqnarray*}
\int_1^{\infty} |\eta_i(y)| \mu(dy) & \leq & \left(\int_1^{\infty}
\, \eta_i^2(y) \mu(dy)\right)^{\frac 12} \, \left(\int_1^{\infty}
\mu(dy)\right)^{\frac 12} < \infty
\end{eqnarray*}
according to Proposition \ref{prop : whyme}. We have thus proved that $\eta_i\in
\L^1(\mu)$.
\end{proof}


We now obtain sharper estimates using  properties of the Dirichlet
heat kernel. For this reason we introduce
$$
\hspace{-3.5cm} \textbf{Hypothesis (H4): }\, \int_{1}^{\infty}e^{-Q(x)}dx<\infty\quad\text{and}\quad
\int_{0}^{1} x\;e^{-Q(x)/2}dx<\infty.
$$

\begin{proposition}\label{propspt}
Assume (H2) and (H4) hold.
Then all eigenfunctions $\eta_k$ belong to $\L^1(\mu)$, and there is a
constant $K_1>0$ such that for any $x\in(0,\infty)$ and any $k$
$$
|\eta_k(x)|\le K_1 \;e^{\lambda_k}\; e^{Q(x)/2}.
$$
Moreover $\eta_1$  is strictly positive on $\R^{+}$, and there is a
constant $K_2>0$ such that for any $x\in(0,1]$ and any $k$
$$
|\eta_{k}(x)|\le K_2\,x   \;  e^{2\lambda_k} \; e^{Q(x)/2}.
$$
\end{proposition}

 \begin{proof}
In Section \ref{sectiongap} we introduced the semigroup $\tilde
P_{t}$ associated with the Schr\"{o}dinger equation and showed
that $\eta_k=e^{\frac Q2}\psi_k$, where $\psi_k$ is the unique
eigenfunction related to  the eigenvalue $\lambda_k$ for $\tilde
P_{t}$. Using estimates on this semigroup, we will get some
properties of $\psi_k$, and we will prove the proposition.

The semigroup $\tilde P_{t}$ is given for
 $f\in \L^{2}\big(\R^{+},dx\big)$   by
$$
\tilde P_{t}f(x)= \, \E^{\W_x} \left[ f(\omega(t)) \,
\BBone_{t<T_0} \, \exp\left( - \, \frac 12 \, \int_0^t \, (q^2 -
q')(\omega_s) ds\right)\right]\;,
$$
where $\E^{\W_x}$ denotes the expectation w.r.t. the Wiener measure
starting from $x$.
 We first
establish a basic estimate on its kernel $\tilde p_{t}(x,y)$.

\begin{lemma}\label{estimkernel}
Assume condition (H2) holds. There exists a constant $K_3>0$ and a
continuous increasing function $B$ defined on $[0,\infty)$
satisfying $\lim_{z\to\infty}B(z)=\infty$, such that for any $x>0,
y>0$ we have
\begin{equation}\label{estim1}
0<\tilde p_{1}(x,y)\le e^{-(x-y)^{2}/4} e^{-B(\max\{x,y\})}\;.
\end{equation}
and
\begin{equation}\label{estim2}
\tilde p_{1}(x,y)\le K_3\; p_{1}^{D}(x,y)\;,
\end{equation}
where $p_{t}^{D}$ is the Dirichlet heat kernel in $\R^{+}$
given for $x,y \in\R^+$ by $$p_{t}^{D}(x,y)={\frac{1}{\sqrt{2\pi
t}}}\left(e^{\frac{-(x-y)^2}{2t}} - e^{-\frac{(x+y)^2}{2t}}\right).$$
\end{lemma}
\medskip
The proof of this lemma is postponed to the Appendix.

\bigskip

\hskip 0.3cm It follows immediately from the previous lemma that
the kernel $\tilde p_{1}(x,y)$ defines a bounded operator $\tilde
P_1$ from $\L^{2}\big(\R^{+},dx\big)$ to
$\L^{\infty}\big(\R^{+},dx\big)$. As a byproduct, we get that all
eigenfunctions $\psi_k$ of $\tilde{P}_1$ are bounded, and more
precisely
$$
|\psi_k|\le K_1 \;e^{\lambda_k}\;.
$$

One also deduces from the previous lemma that the kernel defined
for $M>0$ by
$$
\tilde p_{1}^{M}(x,y)= \BBone_{x<M}\BBone_{y<M}\;\tilde
p_{1}(x,y)\;
$$
is a Hilbert-Schmidt operator in $\L^{2}\big(\R^{+},dx\big)$, in
particular is a compact operator (see for example \cite[pages 177,
267]{conway}). In addition, it follows at once again from Lemma
\ref{estimkernel} that if $\tilde P^{M}_{1}$ denotes the operator
with kernel $\tilde p_{1}^{M}$, we have the following estimate, in
the norm of operators acting on $\L^{2}\big(\R^{+},dx\big)$,
$$
\big \|\tilde P^{M}_{1}- \tilde P_{1}\big
\|_{\L^{2}\big(\R^{+},dx\big)}\le
 C' e^{-B(M)}
$$
where $C'$ is a positive constant independent of $M$. Since
$\lim_{M\to\infty} B(M)=\infty$, the operator  $\tilde P_{1}$ is a
limit in norm of compact operators in $\L^{2}\big(\R^{+},dx\big)$
and hence compact. Since
$\tilde p_{1}(x,y)>0$, the operator $\tilde P_{1}$ is positivity
improving (that is if $0\neq f\ge 0$ then $\tilde P_{1}f>0$ ) implying that  the  eigenvector $\psi_{1}$ is
positive.

  We now claim that $|\psi_k(x)|\leq K_2\:x\;e^{2\lambda_k}$ for
$0<x\le1$. We  have from Lemma
\ref{estimkernel} and the explicit expression for $p^{D}_{1}(x,y)$
the existence of a constant $K_3$ such that
$$
\left|e^{-\lambda_{k}}\psi_{k}(x)\right|\le
K_3\int_{0}^{\infty}p^{D}_{1}(x,y)\,|\psi_{k}(y)|\;dy\le
K_3\|\psi_{k}\|_{\infty} \sqrt{\frac{2}{\pi}} e^{-x^{2}/2}
\int_{0}^{\infty} e^{-y^{2}/2} \sinh(xy)\;dy\;.
$$
We now estimate the integral in the right hand side. Using the
convexity property of $\sinh$ we get  $\sinh(xy)\le x\sinh(y)\le
\frac x2 e^y$, for $x\in[0,1], y\ge 0$ which yields
$$
\int_{0}^{\infty} e^{-y^{2}/2} \sinh(xy)\;dy \le \frac x2
\int_{0}^{\infty} e^{-y^{2}/2} e^y\;dy
$$
proving the claim. Together with hypothesis (H4), this estimate
implies that $\eta_{k}$ belongs to $\L^{1}\big((0,1),d\mu\big)$.
\medskip

Since
$$
\eta_k(x)=\psi_{k}(x)\;e^{Q(x)/2},
$$
we have
$$
\int_{1}^{\infty} \eta_k d\mu= \int_{1}^{\infty}
\psi_{k}(x)\;e^{-Q(x)/2}dx
$$
which implies  $\eta_{k}\in \L^{1}\big((1,\infty),d\mu\big)$
using Cauchy-Schwarz's
 inequality. This finishes the proof of
Proposition \ref{propspt}.
\end{proof}

\begin{remark}
Let us discuss some easy facts about the hypotheses introduced.
\begin{enumerate}
\item[(1)] If $q$ and $q'$ extend continuously up to $0$, hypotheses
(H2), (H3) and (H4) reduce to their counterpart at infinity.
\medskip

\item[(2)] Consider $q(x)=\frac ax + g(x)$ with $g$ a $C^1$
function up to $0$. In order that \eqref{eqexplo2} holds at
the origin we need $a>-{\frac 12}$. Then $\mu(dx) = \Theta(x) x^{-2a} dx$  with
$\Theta$ bounded near
the origin, while $q^2(x)-q'(x)\approx (a+a^2)/x^2$. Hence for
(H2) to hold, we need $a\geq 0$. Now we have the estimates
$$
\int_0^\varepsilon \frac{1}{(q^2(x)-q'(x)+C+2)} \, \mu(dx)\approx
\int_0^\varepsilon \Theta(x) x^{2(1-a)} dx \, ,
$$
and
$$
\int_0^\varepsilon x e^{-{Q(x)\over 2}}dx  \approx
\int_0^\varepsilon \Theta(x) x^{1-a} dx\,.
$$
Therefore (H3) holds for $a<{3\over 2}$ and (H4) holds (at $0$) for $a< 2$. The conclusion is that
$a\in [0,\frac 32)$.
\medskip

\noindent We recall that $a=\frac 12$ if $X$ comes from a
generalized Feller diffusion.

\medskip

\item[(3)] If $q(x)\ge0$ for $x$ large, hypothesis (H2) implies
the first part of hypothesis (H4). Indeed, take $a>0$ be such that for any $x\ge a$ we have
$q(x)>0$ and $q^{2}(x)-q'(x)>1$. Consider the function
$y=e^{-Q/2}$ which satisfies $y'=-qy$ and
$y''=(q^{2}-q')y$. For $b>a$ we get
after integration by parts
$$
\hspace{1cm} 0=\int_{a}^{b} \left((q^{2}-q')y^{2}-yy''\right)\,dx =\int_{a}^{b}
\left((q^{2}-q')y^{2}+{y'}^{2}\right)\,dx -y(b)y'(b)+y(a)y'(a)\;.
$$
Using $y'=-qy$ we obtain
$$
\hspace{0.5cm} \int_{a}^{b} y^2 dx \leq \int_{a}^{b}
\left((q^{2}-q')y^{2}+{y'}^{2}\right)\,dx
=q(a)y(a)^{2}-q(b)y(b)^{2}\leq q(a)y(a)^{2}<\infty
$$
and the result follows by letting $b$ tend to infinity.
\end{enumerate}
\hfill{$\diamondsuit$}
\end{remark}

\bigskip

\section{\bf Quasi-stationary distribution and Yaglom limit}
\label{sectionyaglom}

Existence of the Yaglom limit and of {\it q.s.d.} for killed one-dimensional
diffusion processes have already been proved by various authors,
following the pioneering work by Mandl \cite{Mandl} (see
e.g. \cite{CMM, MSM, SE05} and references therein). One of the main
assumptions in these papers is $\kappa(\infty)=\infty$ and
$$
 \int_1^\infty \, e^{-Q(y)} \,
\left(\int_1^y \, e^{Q(z)} \, dz\right) \, dy \, = \, + \infty
$$
which is not necessarily satisfied in our case. Indeed, under mild
conditions, the Laplace method yields that $\int_1^y \, e^{Q(z)} \, dz$
behaves like $e^{Q(y)}/2q(y)$ when $y$ tends to infinity,
 so the above equality will
not hold if $q$ grows too fast to infinity at infinity. Actually, we
will be particularly interested in these cases (our forthcoming
assumption (H5)), since they are exactly those when the diffusion
``comes down from infinity'', which ensures uniqueness of the {\it q.s.d.}. The
second assumption in the aforementioned papers
 is that $q$ is $C^1$ up to the origin which
is not true in our case of interest.

It is useful to introduce the following condition
\begin{definition}
 We say that \textbf{Hypothesis (H)} is satisfied if
 (H1) and (H2) hold, and moreover $\eta_1 \in
\L^1(\mu)$ (which is the case for example under (H3) or (H4)).
\end{definition}

We now study the existence of {\it q.s.d.} and Yaglom limit in our framework.
When $\eta_1 \in \L^1(\mu)$, a natural candidate for being a {\it q.s.d.}
is the normalized measure $\eta_1 \mu/\langle \eta_1, 1\rangle_\mu$,
which turns out to be the conditional limit distribution.

\begin{theorem}
\label{thm : Yaglom}
Assume that Hypothesis (H) holds. Then
$$
d\nu_1 = \frac{\eta_1 d\mu
}{\langle \eta_{1},1\rangle_{\mu}}
$$
 is a quasi-stationary distribution, that is
for every $t\geq
0$ and any Borel subset $A$ of $(0,\infty)$,
$$
\P_{\nu_1}(X_t \in A
\, | \, T_0>t) = \nu_1(A) \, .
$$

Also for any $x>0$ and any Borel subset $A$ of $(0,\infty)$,
\begin{eqnarray}
\label{limit} \lim_{t \to \infty}e^{\lambda_1t} \,
\P_x( T_0 > t)\, = \, \eta_1(x) \,\langle \eta_{1},1\rangle_{\mu}\,
,\end{eqnarray}
$$
\lim_{t \to \infty} e^{\lambda_1t}\, \P_x(X_t \in A \, , \, T_0 > t)
\, = \, \nu_1(A)\,\eta_1(x) \langle \eta_{1},1\rangle_{\mu}\,.
$$
This implies since $\eta_{1}>0$ on $(0,\infty)$
$$
\lim_{t \to \infty} \, \P_x(X_t \in A \, | \, T_0
> t) \, = \, \nu_1(A) \, ,$$
and the probability measure $\nu_1$ is the Yaglom limit
distribution. Moreover, for any probability measure $\rho$ with compact support
in $(0,\infty)$ we have
\begin{eqnarray}
\label{Yaglom1}
&&\lim_{t \to \infty}e^{\lambda_1t}
\, \P_\rho( T_0 > t) \, =  \langle \eta_{1},1\rangle_{\mu}\int
\eta_1(x) \rho(dx)\;;\\
&&\lim_{t \to \infty} e^{\lambda_1t} \, \P_\rho(X_t \in A \, , \, T_0
> t) \, = \, \nu_1(A)\, \langle \eta_{1},1\rangle_{\mu}\,
\int \eta_1(x)\, \rho(dx)
\;;\\
&&\lim_{t \to \infty} \P_\rho(X_t \in A  \mid T_0
> t) \, = \, \nu_1(A).\label{unlabel}
\end{eqnarray}
\end{theorem}
\begin{proof}
Thanks to the symmetry of the semigroup, we have for all $f$ in
$\L^2(\mu)$,
$$\int P_t f \eta_1 d\mu = \int f P_t \eta_1 d\mu = e^{-\lambda_1 t}
\, \int f \eta_1 d\mu \, .$$ Since $\eta_1 \in \L^1(\mu)$, this
equality extends to all bounded $f$.  In particular we may use it with
$f=\BBone_{(0,\infty)}$  and with
$f=\BBone_A$ . Noticing that $$\int P_t
(\BBone_{(0,\infty)}) \, \eta_1 d\mu = \P_{\nu_1}(T_0>t)\langle
\eta_1, 1\rangle_\mu$$ and $\int P_t \BBone_A \eta_1 d\mu = \P_{\nu_1}(X_t
\in A \, , \, T_0>t)\langle \eta_1, 1\rangle_\mu $, we have shown that
$\nu_1$ is a {\it q.s.d.}.

The rest of the proof is divided into two cases. First assume that
$\mu$ is a bounded measure.
Thanks to Theorem \ref{thmloit}, we know that for any $x>0$, any set
$A \subset (0,\infty)$ such
that $\BBone_A \in \L^2(\mu)$ and for any $t>1$
\begin{eqnarray*}
\P_x(X_t \in A \, , \, T_0>t) & = & \int \, \P_y(X_{t-1} \in A \,
, \, T_0>t-1) \, r(1,x,y) \, \mu(dy)
\\ & = & \int P_{t-1} (\BBone_A)(y) \, r(1,x,y) \, \mu(dy) \\
& = & \int \, \BBone_A(y) \, (P_{t-1}r(1,x,.))(y) \, \mu(dy) \, .
\end{eqnarray*}
Since both $\BBone_A$ and $r(1,x,.)$ are in $\L^2(\mu)$ and since
(H2) is satisfied, we obtain using Theorem \ref{thmspectrum}
\begin{equation}\label{eqyaglom}
\lim_{t \to \infty} \, e^{\lambda_1 (t-1)} \, \P_x(X_t \in A \, ,
\, T_0>t) = \langle \BBone_A,\eta_1 \rangle_\mu \, \langle
r(1,x,.), \eta_1 \rangle_\mu \, .
\end{equation}
Since
$$
\int r(1,x,y) \eta_1(y) \mu(dy) = (P_1 \eta_1)(x)= e^{-\lambda_1}
\eta_1(x)
$$
we get that $\nu_1$ is the Yaglom limit.
\smallskip

If $\mu$ is not bounded (i.e. $\BBone_{(0,\infty)} \notin
\L^2(\mu)$) we need an additional result to obtain the Yaglom
limit.

\begin{lemma}\label{lemloit}
Assume $\eta_1 \in
\L^1(\mu)$ then for all $x>0$, there exists a locally bounded
function $\Theta(x)$ such that for
all $y>0$ and all $t>1$,
\begin{equation}
\label{equnif}
r(t,x,y) \, \leq \, \Theta(x) \, e^{-\lambda_1 t} \, \eta_1(y) \, .
\end{equation}
\end{lemma}

We postpone the proof of the lemma and indicate how it is used to
conclude the proof of the theorem.

If \eqref{equnif} holds, for $t>1$, $e^{\lambda_1 t} \, r(t,x,.)
\in \L^1(\mu)$ and is dominated by $\Theta(x) \, \eta_1$.
Since $r(1,x,\cdot)\in \L^2(\mu)$ by Theorem \ref{thmloit}, using
Theorem \ref{thmspectrum}  and   writing again
 $r(t,x,.)=P_{t-1}r(1,x,.)$ $\mu$ a.s., we deduce that $\lim_{t
\to \infty} e^{\lambda_1 t} \, r(t,x,.)$ exists in $\L^2(\mu)$
and is equal to
$$
e^{\lambda_1} \langle r(1,x,.), \eta_1 \rangle_\mu \, \eta_1(.) =
\eta_1(x) \, \eta_1(.) \, .
$$

Recall that convergence in $\L^2$ implies almost sure convergence along
subsequences.
Therefore,  for any sequence $t_n \to \infty$ there exists a subsequence
$t'_n$ such that
$$\lim_{n \to \infty} e^{\lambda_1 t'_n} \, r(t'_n,x,y) \, = \,
\eta_1(x) \, \eta_1(y) \, \textrm{ for $\mu$-almost all } y>0 \,
.
$$
Since
$$
\P_x(T_0>t'_n) = \int_0^{\infty} \, r(t'_n,x,y) \mu(dy) \, ,
$$
Lebesgue bounded convergence theorem yields
$$
\lim_{n \to \infty} \, e^{\lambda_1 t'_n} \, \P_x(T_0>t'_n)= \eta_1(x) \,
\int_0^{\infty} \eta_1(y) \mu(dy) \, .
$$
That is \eqref{eqyaglom} holds with $A=(0,\infty)$ for the
sequence $t'_n$. Since the limit does not depend on the
subsequence, $\lim_{t \to \infty} \, e^{\lambda_1 t} \,
\P_x(T_0>t)$ exists and is equal to the previous limit, hence
\eqref{eqyaglom} is still true. The rest of this part follows as
before.

For the last part of the theorem, that is passing from the initial
Dirac measures at every fixed $x>0$ to the compactly supported
case, we just use that $\Theta(\bullet)$ is bounded on compact
sets included in $(0,\infty)$.
\smallskip

{\it{Proof of Lemma \ref{lemloit}.}} According to the
parabolic Harnack's inequality (see for example \cite{Tru}), 
for all $x>0$, one can find
$\Theta_0(x)>0$, which is locally bounded, such that for all
$t>1$, $y>0$ and $z$ with $|z-x|\leq \rho(x) = \frac 12 \wedge
\frac x4$
$$
r(t,x,y) \, \leq \, \Theta_0(x) \, r(t+1,z,y)  \, .
$$
 It follows that
\begin{eqnarray*}
r(t,x,y) & = & \frac{\left(\int_ {|z-x|\leq \rho(x)} r(t,x,y) \eta_1(z)
\mu(dz)\right)}{\left(\int_ {|z-x|\leq \rho(x)} \eta_1(z) \mu(dz)\right) }
\\ & \leq & \Theta_0(x) \, \frac{\left(\int_ {|z-x|\leq \rho(x)} r(t+1,z,y)
\eta_1(z) \mu(dz)\right)}{\left(\int_ {|z-x|\leq \rho(x)}
\eta_1(z) \mu(dz)\right) } \\ & \leq & \Theta_0(x) \, \frac{\left(\int
r(t+1,z,y) \eta_1(z) \mu(dz)\right)}{\left(\int_ {|z-x|\leq
\rho(x)} \eta_1(z) \mu(dz)\right) }\\ & \leq & \Theta_0(x) \, \frac{ e^{-
\lambda_1 (t+1)} \, \eta_1(y)}{\left(\int_ {|z-x|\leq \rho(x)}
\eta_1(z) \mu(dz)\right) } \,,
\end{eqnarray*}
since $P_{t+1}\eta_{1}=e^{-\lambda_{1}(t+1)}\eta_{1}$.
But $\Theta_1(x)=\int_ {|z-x|\leq \rho(x)} \eta_1(z) \mu(dz) > 0$,
otherwise $\eta_1$, which is a solution of the linear o.d.e.
$\frac 12 g'' - q g' + \lambda_1 g =0$ on $(0,\infty)$, would
vanish on the whole interval $|z-x|\leq \rho(x)$, hence on
$(0,\infty)$ according to the uniqueness theorem for linear
o.d.e's. The proof of the lemma is thus completed with
$\Theta=e^{-\lambda_1}\Theta_0/\Theta_1$.
\end{proof}

\medskip
The positive real number
$\lambda_1$ is the natural killing rate of the process. Indeed,
the limit \eqref{limit} obtained in Theorem \ref{thm : Yaglom}
shows that for any $x>0$ and any $t>0$,
$$\lim_{s \to \infty} \, \frac{\P_x(T_0>t+s)}{\P_x(T_0>s)} \, =
\, e^{- \lambda_1 t} \, .$$ Let us also remark that
$$\P_{\nu_1}(T_0>t)=e^{-\lambda_1 t}.$$

\medskip

In order to control the speed of convergence to the Yaglom limit, we
first establish the following lemma.

\begin{lemma}\label{lpl2}
Under conditions (H2) and (H4), the operator $P_1$ is bounded from
$\L^\infty(\mu)$ to $\L^2(\mu)$. Moreover, for any compact
subset $K$ of $(0,\infty)$, there is a constant $C_K$ such that
for any function $f\in\L^1(\mu)$ with support in $K$  we have
$$
\|P_1f\|_{\L^2(\mu)}\le C_K\,\|f\|_{\L^1(\mu)}
$$
\end{lemma}

\begin{proof}
Let $g\in \L^\infty(\mu)$, since
$$
\big|P_1g\big|\le P_1\big|g\big|\le \|g\|_{\L^\infty(\mu)}\;,
$$
we get from (H4)
$$
\int_1^\infty \big|P_1g\big|^2 d\mu \le
\|g\|_{\L^\infty(\mu)}^2 \int_1^\infty  e^{-Q(x)}dx\;.
$$
We now recall that (see Section 3)
$$
P_1g(x)=e^{Q(x)/2}\tilde P_1\left(e^{-Q/2}g\right)(x)\;.
$$
It follows from Lemma \ref{estimkernel} that uniformly in
$x\in(0,1]$ we have (using hypothesis (H4))
$$
\left|\tilde P_1\left(e^{-Q/2}g\right)(x)\right|
\le \mathcal O(1)\;\|g\|_{\L^\infty(\mu)} \int_0^\infty
e^{-Q(y)/2}e^{-y^2/4}y\,dy\le \mathcal O(1)\;\|g\|_{\L^\infty(\mu)} \;.
$$
This implies
$$
\int_0^1 \big|P_1g\big|^2 d\mu
=\int_0^1 \big|\tilde P_1\left(e^{-Q/2}g\right)(x)\big|^2 dx \le
\mathcal O(1)
\|g\|_{\L^\infty(\mu)}^2 \;,
$$
and the first part of the lemma follows. For the second part, we have
from the Gaussian bound of Lemma \ref{estimkernel} that for any $x>0$
and for any $f$ integrable and with support in $K$
$$
\left|\tilde P_1\left(e^{-Q/2}f\right)(x)\right|
\le \mathcal O(1) \int_K e^{-Q(y)/2}e^{-(x-y)^2/2}|f(y)|\,dy
$$
$$
\le \mathcal O(1) \sup_{z\in K} e^{Q(z)/2}
\sup_{z\in K} e^{-(x-z)^2/2}
\int_K  e^{-Q(y)} |f(y)|\,dy
\le \mathcal O(1) \;e^{-x^2/4}
\int_K  e^{-Q(y)} |f(y)|\,dy
$$
since $K$ is compact. This implies
$$
\int_0^\infty \big|P_1f\big|^2 d\mu
=\int_0^\infty \big|\tilde P_1\left(e^{-Q/2}f\right)(x)\big|^2 dx \le
\mathcal O(1)\;
\|f\|^2_{\L^1(\mu)} \;.
$$
\end{proof}

We can now  use the
spectral decomposition of $r(1,x,.)$ to obtain the following
convergence result.

\begin{proposition}\label{ratecv}
Under conditions (H2) and (H4), for all $x>0$ and any measurable
subset $A$ of $(0,\infty)$, we have
\begin{equation}\label{eqrate}
\!\!\!\lim\limits_{t \to \infty} \!e^{(\lambda_2 - \lambda_1)t}
\bigg(\P_x(X_t\! \in \! A \! \mid \! T_0\!>\!t)- \nu_1(A)\bigg)\! =\!
\frac{\eta_2(x)}{\eta_1(x)}\! \left(\!\!\frac{\langle
1,\eta_1\rangle_\mu \langle \BBone_A,\eta_2\rangle_\mu \! -\!
\langle 1,\eta_2\rangle_\mu \langle
\BBone_A,\eta_1\rangle_\mu}{\langle 1,\eta_1\rangle_\mu^2}\!\right)\!\!.
\end{equation}
\end{proposition}

\medskip

\begin{proof}
Let $h$ be a non negative bounded function, with compact support
in $(0,\infty)$. By using the semigroup property, Lemma \ref{lpl2}
and the spectral decomposition for compact  self adjoint
semigroups (see Theorem \ref{thmspectrum}), we have for any $t>2$,
$$
\int \P_{x}(X_t\in A\,,\,T_0>t) h(x) dx=\langle he^Q
\,,\,P_t\BBone_A\rangle_\mu =\langle P_1\big(he^Q\big)
\,,\,P_{(t-2)}P_1\BBone_A\rangle_\mu
$$
$$
=\langle P_1\big(he^Q\big)\,,\,\eta_1\rangle_\mu \,\langle\eta_1\,,\,P_1
\BBone_A\rangle_\mu e^{-\lambda_1(t-2)}
+\langle P_1\big(he^Q\big)\,,\,\eta_2\rangle_\mu \,\langle\eta_2\,,\,P_1
\BBone_A\rangle_\mu e^{-\lambda_2(t-2)}
+R(h,A,t)
$$
with
$$
|R(h,A,t)|\le\sum\limits_{i\ge 3} e^{-\lambda_i(t-2)}
\Big|\langle P_1\big(he^Q\big),\eta_i\rangle_\mu \langle\eta_i,P_1
\BBone_A\rangle_\mu \Big|\le
e^{-\lambda_3(t-2)}\|P_1\big(he^Q\big)\|_{\L^2(\mu)}
\|P_1\BBone_A\|_{\L^2(\mu)}\;,
$$
due to $\lambda_{1}< \lambda_{2}< \lambda_{3}< \ldots$,
the Cauchy-Schwarz inequality
and Parseval's identity. Note that since $P_1$ is symmetric with
respect to the scalar product, we have $\langle
P_1\big(he^Q\big)\,,\,\eta_1\rangle_\mu=e^{-\lambda_1}\langle
he^Q,\eta_1\rangle_\mu$ and similarly  for $\eta_2$. We also have
$\langle\eta_1\,,\,P_1 \BBone_A\rangle_\mu =e^{-\lambda_1}
\langle\eta_1\,,\,\BBone_A\rangle_\mu$
  and similarly for $\eta_2$.
It follows immediately from Lemma \ref{lpl2} that for any fixed compact
subset $K$ of $(0,\infty)$, any $A$ and any $h$ satisfying
the hypothesis of the proposition with support contained in $K$,
$$
|R(h,A,t)|\le \mathcal O(1) e^{-\lambda_3(t-2)}\|he^Q\|_{\L^1(\mu)}
\le \mathcal O(1) e^{-\lambda_3(t-2)}\|h\|_{\L^1(dx)}
\;,
$$
since $h$ has compact support in $(0,\infty)$.
Therefore, letting $h$ tend to a Dirac mass, we obtain that for  any
compact subset $K$ of $(0,\infty)$, there is a constant $D_K$  such that
for any $x\in K$, for any measurable subset $A$ of $(0,\infty)$,
 and for any $t>2$, we have
$$
\left|\P_x(X_t\in A\,,\,T_0>t)-
e^{Q(x)}
\eta_1(x)\,\langle\eta_1\,,\,\BBone_A\rangle_\mu e^{-\lambda_1t}
-e^{Q(x)}
\eta_2(x)\,\langle\eta_2\,,\,\BBone_A\rangle_\mu e^{-\lambda_2t}\right|
\le D_K\; e^{-\lambda_3t}\;.
$$
The proposition follows at once from
$$
\P_x(X_t\in A\,|\,T_0>t)=\frac{\P_x(X_t\in A\,,\,T_0>t)}{
\P_x(X_t\in (0,\infty)\,,\,T_0>t)}\;.
$$
\end{proof}

\section{\bf The $Q$-process}\label{secQ}

 As in \cite{CMM} (Theorem B), we can also describe the
law of the process conditioned to be never  extinct, usually
called the $Q$-process (also see \cite{Lam06}).

\bigskip

\begin{corollary}\label{corq}
Assume  (H) holds.  For all $x>0$ and $s\ge 0$ we have
$$
\lim_{t \to \infty} \, \P_x (X \in B \mid T_0 > t) =
\Q_x(B) \, \hbox{ for all } B  \hbox{ Borel
measurable subsets of } C([0,s]),
$$
where $\Q_x$ is the law of a diffusion process on $(0,\infty)$,
with transition probability densities (w.r.t. Lebesgue measure) given by
$$
q(s,x,y)=
e^{\lambda_1 s} \, \frac{\eta_1(y)}{\eta_1(x)} \, r(s,x,y) \, e^{-Q(y)} \, ,
$$
that is, $\Q_x$ is locally absolutely continuous w.r.t. $\PP_x$ and
$$
\Q_x(X\in B) = \E_x\left(\BBone_{B}(X) \,
e^{\lambda_1 s}\, \frac{\eta_1(X_s)}{\eta_1(x)}\, , T_0>s\right)\, .
$$
\end{corollary}
\begin{proof}
First check thanks to Fubini's theorem and $\kappa(0^+)<\infty$ in
Hypothesis (H1), that $\Lambda(0^+)> - \infty$. We can thus slightly
change the notation (for this proof only) and define $\Lambda$ as
$\Lambda(x)=\int_0^x e^{Q(y)}\, dy$. From standard diffusion theory,
$(\Lambda(X_{t\wedge T_0}) ; t\ge 0)$ is a local martingale, from which it is easy
to derive that for any $y\ge x \ge 0$, $\PP_y (T_x <T_0) =
\Lambda(y)/\Lambda(x).$

Now define $v(t,x)=\frac{\P_x(T_0>t)}{\P_1(T_0>t)}$. As in \cite[proof
of Theorem B]{CMM}, one can prove for any $x\ge 1$, using the strong
Markov property at $T_x$ of the diffusion $X$ starting from 1, that
$v(t,x)\le \Lambda(x)/\Lambda(1)$. On the other hand for $x\le 1$ we obtain
$v(t,x)\le 1$. Then  for any $x\ge 0$ we get $v(t,x)\le 1+\Lambda
(x)/\Lambda(1)$.

Now thanks to Theorem \ref{thm : Yaglom}, for all $x$,
$e^{\lambda_1 t} \P_x(T_0>t) \to
 \eta_1(x)\langle 1,\eta_{1}\rangle_{\mu}$ as $t \to \infty$, and

$$
\lim_{t \to \infty} v(t,x)  = \frac{\eta_1(x)}{\eta_1(1)} \, .
$$
Using the Markov property, it is easily seen that for large $t$,
$$\P_x ( X \in B \mid T_0 > t) = \E_x\left[\BBone_{B}(X) \,
v(t-s,X_s), T_0>s\right] \, \, \frac{\P_1(T_0> t-s)}{\P_x(T_0>t)} \, .$$ The random variable in
the expectation is (positive and) bounded from above by $1+\Lambda(X_s)/\Lambda(1)$, which is
integrable (see below), so we obtain the desired result using Lebesgue bounded convergence
theorem.

To see that $\E_x\left(\Lambda(X_s) \, \BBone_{s<T_0}\right)$ is
finite, it is enough to use It\^o's formula with the harmonic
function $\Lambda$ up to time $T_0\wedge T_M$. Since $\Lambda$ is
nonnegative it easily yields $\E_x\left(\Lambda(X_s) \,
\BBone_{s<T_0\wedge T_M}\right) \leq \Lambda(x)$ for all $M>0$.
Letting $M$ go to infinity the indicator converges almost surely to $\BBone_{s<T_0}$ 
(thanks to Hypothesis (H1)) so  the monotone convergence
theorem yields $\E_x\left(\Lambda(X_s) \, \BBone_{s<T_0}\right)
\leq \Lambda(x)$.
\end{proof}

Recall that $\nu_1$ is the Yaglom limit.
\begin{corollary}\label{corjavaitort}
Assume (H) holds. Then for any Borel subset $B\subseteq (0,\infty)$ and any $x$,
$$\lim_{s \to \infty} \, \Q_x(X_s \in B) = \int_B \, \eta^2_1(y) \mu(dy)\,=
\langle\eta_{1},1\rangle_{\mu} \int_B \, \eta_1(y)  \nu_1(dy) .
$$

\end{corollary}

\begin{proof}
We know from the proof of Theorem \ref{thm : Yaglom}
 that $e^{\lambda_1 s} r(s,x,.)$ converges to $\eta_1(x) \,
\eta_1(.)$ in $\L^2(\mu)$ as $s\to\infty$.
Hence, since $\BBone_B \eta_1 \in
\L^2(\mu)$,
$$
\eta_1(x)\Q_x(X_s\in B)=\int \BBone_B(y) \eta_1(y) \,
e^{\lambda_1 s} r(s,x,y) \mu(dy) \to \eta_1(x) \, \int_B
\eta^2_1(y) \mu(dy)
$$
as $s \to \infty$.
We remind the reader  that
$d\nu_{1}=\eta_{1}d\mu/\langle\eta_{1},1\rangle_{\mu}$.
\end{proof}
\begin{remark}
The stationary
measure of the $Q$-process is absolutely continuous w.r.t. $\nu_1$,
with Radon-Nikodym derivative $\langle\eta_{1},1\rangle_{\mu}\eta_1$
 which thanks to Proposition
\ref{prop : whyme} is \emph{nondecreasing}. In particular, the
ergodic measure of the $Q$-process dominates stochastically the Yaglom
limit.  We refer to \cite{MSM,Lam06} for further discussion of the
relationship between {\it q.s.d.} and ergodic measure of the $Q$-process.
\hfill{$\diamondsuit$}
\end{remark}

\section{\bf Infinity is an entrance boundary and uniqueness of {\it q.s.d.}}
\label{sectionattrac}

We start with the notion of quasi limiting distribution {\it q.l.d.}.

\begin{definition} A probability measure $\pi$ supported on $(0,\infty)$ is a {\it q.l.d.} if there exists a probability
measure $\nu$ such that the following limit exists in distribution
$$
\lim\limits_{t\to \infty} \P_\nu( X_t \in \bullet \mid T_0>t)=\pi(\bullet).
$$
We also say that $\nu$ is attracted to $\pi$, or is in the domain of attraction of $\pi$,
for the conditional evolution.
\end{definition}

Obviously every {\it q.s.d.} is a {\it q.l.d.}, because such measures are fixed points
for the conditional evolution.  We prove that the reciprocal is also true, so both concepts
coincide.

\begin{lemma}
\label{qlid-qsd} Let $\pi$ a probability measure supported on $(0,\infty)$. If $\pi$ is a
{\it q.l.d.} then $\pi$ is a {\it q.s.d.}.
In particular there exists $\alpha \ge 0$ such that for all $s>0$
$$
\P_\pi(T_0>s)=e^{-\alpha s}.
$$
\end{lemma}

\begin{proof} By hypothesis there exists a probability measure $\nu$ such that
$\lim_{t\to \infty} \P_\nu(X_t \in \bullet \mid T_0>t)=\pi(\bullet)$, in distribution. That is, for all
continuous and bounded functions $f$ we have
$$
\lim\limits_{t\to \infty} \frac{\P_\nu(f(X_t), T_0>t)}{\P_\nu( T_0>t)}= \int f(x) \pi(dx).
$$
If we take $f(x)=\P_x( X_s \in A, T_0>s)$, since $f(x)=\int_A r(t,x,y) \mu(dy)$, an application of Harnack's inequality and of the dominated convergence theorem ensures that $f$ is continuous in $(0,\infty)$. 

First, take $A=(0,\infty)$, so that $f(x)=\P_x(T_0>s)$. Then, we obtain for all $s\ge 0$
$$
\lim\limits_{t \to \infty} \frac{\P_\nu(T_0>t+s)}{\P_\nu( T_0>t)}= \P_\pi(T_0>s).
$$
The left hand side is easily seen to be exponential in $s$ and then
there exists $\alpha\ge 0$ such that
$$
\P_\pi(T_0>s)=e^{-\alpha s}.
$$
Second, take $f(x)=\P_x( X_s \in A, T_0>s)$ to conclude that
$$
\begin{array}{ll}
\P_\pi(X_s \in A, T_0>s)&= \lim\limits_{t\to \infty} \P_\nu( f(X_t) | T_0>t) =
\lim\limits_{t\to \infty} \P_\nu( X_{t+s}\in A \mid T_0>t+s) \frac{\P_\nu(T_0>t+s)}{\P_\nu(T_0>t)}\\
&=e^{-\alpha s} \pi(A),
\end{array}
$$
and then $\pi$ is a {\it q.s.d.}.
\end{proof}

Recall from Theorem \ref{thm : Yaglom} that under Hypothesis (H),
the measure $d\nu_1 = \eta_1 d\mu / \langle\eta_{1},1\rangle_{\mu}$ is the Yaglom limit,
which in addition is a {\it q.l.d.} attracting all initial distribution with compact
support on $(0,\infty)$.
It is natural to ask about the uniqueness of the {\it q.s.d.}.
Here again, our assumptions on the behavior of $q$
at infinity will allow us to characterize the domain of attraction
of the {\it q.s.d.} $\nu_1$  associated to $\eta_1$. This turns out to be entirely
different from  the cases studied in \cite{CMM} for instance.

\bigskip

We say that the diffusion process $X$ \emph{comes down from infinity}
if  there is $y>0$ and a time $t>0$ such that
$$
\lim_{x\uparrow \infty} \P_x(T_y<t)\,>\, 0.
$$
This terminology is equivalent to the property that $\infty$ is an entrance boundary for $X$ (for instance see \cite[page 283]{RY}).
\medskip

Let us introduce the following condition
$$
\hspace{-5.5cm} \textbf{Hypothesis (H5): } \, \int_1^\infty \, e^{Q(y)} \,
\int_y^\infty \, e^{-Q(z)} \, dz \, dy \, < \,  \infty.
$$
Tonelli's Theorem ensures that (H5) is equivalent to
\begin{equation}
\label{entrance}
\int_1^\infty  e^{-Q(y)} \int_1^y e^{Q(z)} dz\; dy<\infty.
\end{equation}
\medskip
If (H5) holds
then for $y\ge 1$, $\int_y^\infty \, e^{-Q(z)} \, dz<\infty$. Applying the Cauchy-Schwarz
inequality we get
$(x-1)^2=\left(\int_1^x e^{Q/2}e^{-Q/2} dz\right)^2 \le \int_1^x e^{Q} dz \; \int_1^x e^{-Q} dz$,
and therefore (H5) implies that $\Lambda(\infty)=\infty$.

\medskip

Now we state the main result of this section.

\begin{theorem}
\label{thm : uniqueness}
 Assume  (H) holds.
 Then the following are equivalent:
 \begin{itemize}
\item[(i)] $X$ comes down from infinity;
\item[(ii)] (H5) holds;
\item[(iii)] $\nu_1$ attracts all initial distributions $\nu$ supported in $(0,\infty)$, that is
$$
\lim\limits_{t \to \infty} \PP_\nu(X_t \in \bullet \mid T_0>t)=\nu_1 (\bullet).
$$
\end{itemize}
In particular any of these three conditions implies that there is a unique {\it q.s.d.}.

\end{theorem}

\begin{remark}\label{rempatrick}
It is not obvious when Condition (H5) holds.  In this direction, the following explicit
conditions on $q$, all together, are sufficient for (H5) to hold:
\begin{itemize}
\item $q(x)\geq q_0 > 0$ for all $x \geq x_0$
\item $\limsup_{x\to \infty} q'(x)/2 q^2(x)<1$
\item $\int_{x_0}^{\infty} \,
\frac{1}{q(x)} \, dx < \infty \, .$
\end{itemize}
Indeed, check first that these conditions imply that $q(x)$ goes to
infinity as $x\to\infty$. Then defining  $s(y):= \int_y^\infty e^{-Q(z)}\,
dz$, the first condition above implies  that $s(y)e^{Q(y)}$ is
bounded in $y\ge x_{0}$. Integrating by parts on $\int se^Q dz$ gives
$$
\int_{x_0}^x se^Q dz=\int_{x_0}^x \frac{s}{2q} \; e^Q 2q\, dz=\frac{s}{2q} e^Q\Big|_{x_0}^x
+\int_{x_0}^x \frac{1}{2q} dz + \int_{x_0}^x se^Q \frac{q'}{2q^2} dz.
$$
Since $se^Q/2q$
vanishes at infinity, the third condition
implies that $se^Q(1-q'/2q^2)$ is integrable and  thanks to the
second condition we conclude that (H5) holds.

On the other hand, if (H5) holds, $q'(x)\ge 0$
for $x\ge x_0$ and $q(x_0)>0$, then $q(x)$ goes
to infinity as $x\to\infty$ and $\int_{x_0}^{\infty} \,
\frac{1}{q(x)} \, dx < \infty \, .$

We can retain that under the assumption that $q'(x)\ge 0$ for $x\ge x_0$
and $q(x)$ goes to infinity as $x\to\infty$, then
$$
\mbox{(H5)}
\Longleftrightarrow\int_1^{\infty} \,
\frac{1}{q(x)} \, dx < \infty \, .
$$
Indeed, the only thing left to prove is the sufficiency of (H5). Since $s(y)$
tends to $0$ as $y\to \infty$ (because $Q$ grows at least linearly),
then by the mean value theorem
we have
$$
\frac{\int_y^\infty e^{-Q(z)} dz}{e^{-Q(y)}}=\frac{1}{2q(\xi)},
$$
where $\xi\in [y,\infty)$. Using that $q$ is monotone we obtain the bound
$$
\frac{\int_y^\infty e^{-Q(z)} dz}{e^{-Q(y)}}\le \frac{1}{2q(y)},
$$
and the equivalence is shown.
\hfill{$\diamondsuit$}
\end{remark}

\bigskip

The proof of Theorem \ref{thm : uniqueness} follows  from Propositions  \ref{procmdfi},
\ref{proexp} and \ref{pro : pas oim}.

\begin{proposition}
\label{procmdfi} Assume (H1) holds. If there is a unique {\it q.s.d.} that
attracts all initial distributions supported in $(0,\infty)$,
then $X$ comes down from infinity.
\end{proposition}

\begin{proof} Let $\pi$ be the unique {\it q.s.d.} that attracts all distributions.
We know that  $\P_{\pi}(T_0>t)= e^{-\alpha t}$ for some $\alpha\ge 0$. Since
absorption is certain then $\alpha>0$.  For the rest of the
proof let $\nu$ be any initial distribution supported on $(0,\infty)$, which by
hypothesis is in the domain
of attraction of $\pi$ that is, for any bounded and continuous function $f$ we have
$$
\lim_{t\to\infty}\int_{0}^\infty \P_\nu ( X_t \in dx \mid T_0>t) f(x) = \int_0^\infty  f(x)
\pi(dx).
$$
We now prove that for any $\lambda<\alpha$,
$\E_\nu(e^{\lambda T_0})<\infty$.
As in Lemma \ref{qlid-qsd} we have for any $s$
$$
\lim_{t\to\infty}\frac{ \P_\nu ( T_0>t+s)}{\P_\nu ( T_0>t)} = e^{-\alpha s}\, .
$$ Now pick $\lambda\in(0,\alpha)$ and $\varepsilon>0$ such that
$(1+\varepsilon)e^{\lambda-\alpha}< 1$. An elementary induction shows
that there is $t_0$ such that for any $t>t_0$, and any integer $n$
$$
\frac{ \P_\nu ( T_0>t+n)}{\P_\nu ( T_0>t)}\le
(1+\varepsilon)^{n}e^{-\alpha n}\, .
$$
Breaking down the integral $\int_{t_0}^\infty \P_\nu(T_0>s) \,
e^{\lambda s}\, ds$ over intervals of the form $(n, n+1]$ and using
the previous inequality, it is easily seen that this integral
converges. This proves that $\E_\nu(e^{\lambda T_0})<\infty$ for any
initial distribution $\nu$.

\smallskip

Now fix $\lambda = \alpha/2$ and for any $x\ge 0$, let $g(x)=\E_x(e^{\lambda T_0})<\infty$. We
want to show that $g$ is bounded, which trivially entails that $X$ comes down from infinity.
Thanks to the previous step, for any nonnegative random variable $Y$ with law $\nu$
$$
\E (g(Y)) = \E_\nu (e^{\lambda T_0}) <\infty.
$$
Since $Y$ can be any random variable, this implies that $g$ is bounded. Indeed, observe that $g$
is increasing and $g(0)=1$, so that $a:=1/g(\infty)$ is well defined in $[0,1)$. Then check that
$$
\nu(dx) = \frac{g'(x)}{(1-a)g(x)^2} \, dx
$$
is a probability density on $(0,\infty)$. To conclude we use the fact that $\int g \, d\nu<\infty$ to get
$$
\int g \, d\nu=\int_0^\infty \frac{g'(x)}{(1-a)g(x)} \, dx=\frac{1}{1-a}\ln g(x)\big|_0^\infty=\frac{\ln g(\infty)}{1-a},
$$
and then $g$ is  bounded.

\end{proof}

\begin{proposition}\label{proexp}
The following are equivalent
\begin{itemize}
\item[(i)] $X$ comes down from infinity;
\item[(ii)] (H5) holds;
\item[(iii)] for any $a>0$ there exists $y_a>0$ such that $\sup_{x>y_a} \E_x[e^{aT_{y_a}}] < \infty$.
\end{itemize}
\end{proposition}

\begin{proof} Since $(i)$ is equivalent to $\infty$ being an entrance boundary and
$(ii)$ is equivalent to (\ref{entrance}) we must show
that ``$\infty$ is an entrance boundary'' and (\ref{entrance}) are equivalent.
This will follow from \cite[Theorem 20.12,$(iii)$]{kallenberg}. For that purpose
consider $Y_t=\Lambda(X_t)$. Under each one of the conditions $(i)$ or $(ii)$ we have $\Lambda(\infty)=\infty$.
It is direct to prove that $Y$ is in natural scale on the interval
$(\Lambda(0),\infty)$, that is, for $\Lambda(0)<a\le y \le b <\infty=\Lambda(\infty)$
$$
\P_y(T^Y_a<T^Y_b)=\frac{b-y}{b-a},
$$
where $T^Y_a$ is the hitting time of $a$ for the diffusion $Y$.
Then, $\infty$ is an entrance boundary for $Y$ if and only if
$$
\int_0^\infty y \; m(dy)<\infty,
$$
where $m$ is the speed measure of $Y$, which is given by
$$
m(dy)=\frac{2\;dy }{(\Lambda'(\Lambda^{-1}(y))^2},
$$
see \cite[formula (5.51)]{karatzas}, because $Y$ satisfies the SDE $$
dY_t=\Lambda'(\Lambda^{-1}(Y_t)) dB_t.
$$
After a change of variables we obtain
$$
\int_0^\infty y \; m(dy)=\int_1^\infty e^{-Q(y)} \int_1^x e^{Q(z)} dz\; dx.
$$
Therefore we have shown the equivalence between $(i)$ and $(ii)$.

We continue the proof with $(ii)\Rightarrow (iii)$.
Let $a>0$, and pick $x_a$ large enough so that
$$
\int_{x_a}^\infty e^{Q(x)} \int_x^\infty \,e^{-Q(z)} \, dz\, dx
\le \frac{1}{2a}\, .
$$
Let $J$ be the nonnegative increasing function defined on
$[x_a,\infty)$ by
$$
J(x)=\int_{x_a}^x e^{Q(y)}\int_y^\infty \, e^{-Q(z)} \, dz\, dy.
$$
Then check that $J''= 2q J' -1$, so that $LJ = -1/2$. Set now
$y_a=1+x_a$, and consider a large $M>x$. It\^o's formula gives
\begin{eqnarray*}
\E_x(e^{a(t\wedge T_M\wedge T_{y_a})} \, J(X_{t\wedge T_M\wedge
T_{y_a}})) = J(x) + \E_x \left(\int_0^{t\wedge T_M\wedge T_{y_a}}
\, e^{as}\,(a J(X_s)+LJ(X_s)) \, ds\right) \, .
\end{eqnarray*}
But $LJ=-1/2$, and $J(X_s)< J(\infty) \le 1/(2a)$ for any $s\leq
T_{y_a}$, so that
\begin{eqnarray*}
\E_x[e^{a(t\wedge T_M\wedge T_{y_a})} \, J(X_{t\wedge T_M\wedge
T_{y_a}})] \leq J(x) \,  .
\end{eqnarray*}
But $J$ is increasing, hence for $x\geq y_a$ one gets $1/(2a)>J(x)\geq
J(y_a)> 0$. It follows that $\E_x(e^{a(t\wedge T_M\wedge
T_{y_a})}) \leq 1/(2a J(y_a))$ and finally $\E_x(e^{a T_{y_a}})
\leq 1/(2a J(y_a))$, by the monotone convergence theorem. So $(iii)$ holds.

Finally, it is clear that $(iii) \Rightarrow (i)$.
\end{proof}

\begin{proposition}
\label{pro : pas oim}
Assume  (H) holds. If there is $x_0$ such that $\sup_{x\ge x_0}\E_x(e^{\lambda_1 T_{x_0}})<\infty$,
then $\nu_1$ attracts all initial distribution supported in $(0,\infty)$.
\end{proposition}

The proof of this result requires the following control near $0$ and $\infty$.

\begin{lemma}
\label{tension} Assume (H) holds, and
$\sup_{x\ge x_0}\E_x(e^{\lambda_1 T_{x_0}})<\infty$.
For $h\in \L^1(\mu)$ strictly positive in $(0,\infty)$
we have
\begin{eqnarray}
\label{tight1}
&\lim\limits_{\epsilon\downarrow 0}\limsup\limits_{t\to \infty}
\frac{\int_0^\epsilon h(x) \PP_x(T_0>t) \mu(dx)}
{\int h(x) \PP_x(T_0>t) \mu(dx)}=0\\
\label{tight2}
&\lim\limits_{M\uparrow \infty}\limsup\limits_{t\to \infty}
\frac{\int_M^\infty h(x) \PP_x(T_0>t) \mu(dx)}
{\int h(x) \PP_x(T_0>t) \mu(dx)}=0
\end{eqnarray}
\end{lemma}

\begin{proof}
We start with \eqref{tight1}. Using Harnack's inequality,
we have for $\epsilon<1$ and large $t$

$$
\frac{\int_0^\epsilon h(x) \PP_x(T_0>t) \mu(dx)}
{\int h(x) \PP_x(T_0>t) \mu(dx)}\le \frac{\PP_1(T_0>t)\int_0^\epsilon h(z) \mu(dz)}
{C r(t-1,1,1) \int_1^2 h(x) \mu(dx) \int_1^2 \mu(dy)},
$$
then
$$
\limsup\limits_{t\to \infty}
\frac{\int_0^\epsilon h(x) \PP_x(T_0>t) \mu(dx)}{\int h(x) \PP_x(T_0>t) \mu(dx)}\le
\limsup\limits_{t\to \infty} \frac{\PP_1(T_0>t)\int_0^\epsilon h(z) \mu(dz)}
{C r(t-1,1,1) \int_1^2 h(x) \mu(dx) \int_1^2 \mu(dy)}
$$
$$
=
\frac{e^{-\lambda_1} \langle\eta_{1},1\rangle_{\mu} \int_0^\epsilon h(z) \mu(dz)}
{C\;\eta_1(1)\int_1^2 h(x) \mu(dx) \int_1^2 \mu(dy)},
$$
and the first assertion of the lemma is proven.

For the second limit, we set
$A_{0}:=\sup\limits_{x\ge x_0} \EE_x(e^{\lambda_1 T_{x_0}})<\infty$.
Then for large $M>x_0$, we have
$$
\PP_x(T_0>t)=\int_0^t \PP_{x_0} (T_0>u) \PP_x(T_{x_0} \in d(t-u))
+ \PP_x(T_{x_0}>t).
$$
Using that $\lim\limits_{u\to \infty} e^{\lambda_1 u}
\PP_{x_0}( T_0>u)=\eta_1(x_0)\langle\eta_{1},1\rangle_{\mu}$ we obtain that
$B_{0}:=\sup\limits_{u\ge 0}e^{\lambda_1 u} \PP_{x_0}( T_0>u)<\infty$. Then
$$
\begin{array}{ll}
\PP_x(T_0>t)&\le B_{0}\int_0^t e^{-\lambda_1 u} \PP_x(T_{x_0} \in
  d(t-u))  + \PP_x(T_{x_0}>t)\\
&\le B_{0} e^{-\lambda_1 t} \EE_x(e^{\lambda_1 T_{x_0}})+
e^{-\lambda_1 t}\EE_x(e^{\lambda_1 T_{x_0}})
\le e^{-\lambda_1 t}A_{0}(B_{0}+1),
\end{array}
$$
and \eqref{tight2} follows immediately.
\end{proof}

\bigskip

\begin{proof}[Proof of Proposition \ref{pro : pas oim}] Let $\nu$ be any fixed
probability distribution whose support is contained in $(0,\infty)$. We must show
that the conditional evolution of $\nu$ converges to $\nu_1$.
We begin by claiming that  $\nu$ can be assumed to have a
strictly positive density $h$, with respect to $\mu$.
Indeed, let
$$
\ell(y)=\int r(1,x,y) \nu(dx).
$$
Using Tonelli's theorem we have
$$
\int\int r(1,x,y) \nu(dx)\, \mu(dy)=\int\int r(1,x,y) \, \mu(dy)\,
\nu(dx)=\int \PP_x(T_0>1) \nu(dx)\le 1,
$$
which implies that $\int r(1,x,y) \nu(dx)$ is finite $dy-$a.s.. Also
$\ell$ is strictly positive by  Harnack's inequality. Finally, define
$h=\ell/\int \ell d\mu$. Notice that for $d\rho=h d\mu$
$$
\P_\nu(X_{t+1} \in \bullet \mid T_0>t+1)=\P_\rho(X_t \in \bullet \mid T_0>t),
$$
showing the claim.

Consider $M>\epsilon>0$ and  any
Borel set $A$ included in $(0,\infty)$. Then
$$
\left|\frac{\int \PP_x(X_t \in A,\, T_0>t)  h(x)\,\mu(dx)}{\int
 \PP_x( T_0>t)  h(x)\,\mu(dx)}-
\frac{\int_\epsilon^M  \PP_x(X_t \in A,\, T_0>t) h(x)\,
  \mu(dx)}{\int_\epsilon^M  \PP_x( T_0>t)  h(x)\,\mu(dx)}\right|
$$
is bounded by the sum of the following two terms
$$
\begin{array}{ll}
I1&=\left|\frac{\int \PP_x(X_t \in A,\, T_0>t) h(x)\,  \mu(dx)}{\int
   \PP_x( T_0>t) h(x)\, \mu(dx)}-
\frac{\int_\epsilon^M \PP_x(X_t \in A,\, T_0>t) h(x)\, \mu(dx)}{\int
 \PP_x( T_0>t)  h(x)\,\mu(dx)}\right|\\
&\\
I2&=\left|\frac{\int_\epsilon^M \PP_x(X_t \in A,\, T_0>t)
  h(x)\,\mu(dx)}{\int \PP_x( T_0>t)  h(x)\,\mu(dx)}-
\frac{\int_\epsilon^M  \PP_x(X_t \in A,\, T_0>t)  h(x)\,
\mu(dx)}{\int_\epsilon^M  \PP_x( T_0>t)  h(x)\,\mu(dx)}\right|.
\end{array}
$$
We have the bound
$$
I1\vee I2\le \frac{\int_0^\epsilon  \PP_x( T_0>t) h(x)\,
 \mu(dx)+\int_M^\infty  \PP_x( T_0>t)  h(x)\,\mu(dx)}
{\int \PP_x( T_0>t)  h(x)\,\mu(dx)}.
$$
Thus, from Lemma \ref{tension} we get
$$
\lim\limits_{\epsilon\downarrow 0, \, M\uparrow \infty}
\limsup_{t\to \infty} \left|\frac{\int  \PP_x(X_t \in A,\, T_0>t)
  h(x)\,
\mu(dx)}{\int  \PP_x( T_0>t)  h(x)\,\mu(dx)}-
\frac{\int_\epsilon^M  \PP_x(X_t \in A,\, T_0>t) h(x)\,
 \mu(dx)}{\int_\epsilon^M  \PP_x( T_0>t)  h(x)\,\mu(dx)}\right|=0.
$$
On the other hand we have using (\ref{unlabel})
$$
\lim\limits_{t\to \infty} \frac{\int_\epsilon^M  \PP_x(X_t \in A,\,
  T_0>t)  h(x)\, \mu(dx)}
{\int_\epsilon^M \PP_x( T_0>t)  h(x)\,
\mu(dx)}=\frac{\int_A \eta_1(z) \mu(dz)}{\int_{\RR^+} \eta_1(z)
  \mu(dz)}=\nu_{1}(A),
$$
independently of $M>\epsilon>0$, and the result follows.
\end{proof}

The following corollary of Proposition \ref{proexp} describes how fast the
process comes down from infinity.
\begin{corollary}\label{corexp} Assume (H) and (H5) hold.
Then for all $\lambda < \lambda_1$, $\sup_{x>0}
\E_x[e^{\lambda T_0}] < \infty$.\end{corollary}

\begin{proof} We have seen in  Section \ref{sectionyaglom} (Theorem
  \ref{thm : Yaglom})
that for all $x>0$, $\lim_{t\to \infty} e^{\lambda_1 t}
\P_x\big(T_0>t\big) = \eta_1(x)\langle\eta_{1},1\rangle_{\mu}
<\infty$ i.e. $\E_x[e^{\lambda T_0}]<
\infty$ for all $\lambda < \lambda_1$. Applying Proposition \ref{proexp}
with $a=\lambda$ and the strong Markov property it follows that
$\sup_{x>y_\lambda} \E_x[e^{\lambda T_0}] < \infty$. Furthermore,
thanks to the uniqueness of the solution of \eqref{eqdiff}, $X_t^x
\le  X_t^{y_\lambda}$ a.s. for all $t>0$ and all $x<y_\lambda$, hence
$\E_x[e^{\lambda T_0}]\leq \E_{y_\lambda}[e^{\lambda T_0}]$ for
those $x$, completing the proof.\end{proof}

\medskip

The previous corollary states that  the  killing time
for the process starting from infinity
has exponential moments up to order $\lambda_1$. In \cite{Lam05}
an explicit calculation of the law of $T_0$ is done in the case of
the logistic Feller diffusion $Z$ (hence the
corresponding $X$) and also for other related models. In
particular it is shown in Corollary 3.10 therein, that the
absorption time for the process starting from infinity has a
finite expectation. As we remarked in studying examples, a very
general family of  diffusion processes (including the
logistic one) satisfy all assumptions in Corollary \ref{corexp},
which is thus an improvement of the quoted result.

\bigskip

We end this section by gathering some known results on birth--death processes that are close
to our findings. Let $Y$ be a birth--death process with birth rate $\lambda_n$ and death rate $\mu_n$
when in state $n$. Assume that $\lambda_0=\mu_0=0$ and that extinction (absorption at 0)
occurs with probability 1.
Let
$$
S=\sum_{i\ge 1}\pi_i+\sum_{n\ge 1}(\lambda_n \pi_n)^{-1}\sum_{i\ge n+1} \pi_i\, ,
$$
where
$$
\pi_n=\frac{\lambda_1\lambda_2\cdots \lambda_{n-1}}{\mu_1\mu_2\cdots\mu_n}\, .
$$
In this context sure absorption at $0$, that is (H1), is equivalent to
$A:=\sum_{i\ge 1} (\lambda_i \pi_i)^{-1}=\infty$ (see \cite[formula 7.9]{karlin}). On the other hand we also have
$\E_1(T_0)=\sum_{i\ge 1} \pi_i$.
We may state

\begin{proposition}
For a birth--death process $Y$ that satisfies  (H1),
the following are equivalent:
\begin{itemize}
\item[(i)] $Y$ comes down from infinity;
\item[(ii)] There is one and only one {\it q.s.d.};
\item[(iii)] $\lim_{n\uparrow\infty}\uparrow \E_n(T_0)<\infty$;
\item[(iv)]$S<\infty$.
\end{itemize}
\end{proposition}
\begin{proof}

In \cite[Theorem 3.2]{VD}, it is stated the following key alternative: $S<\infty$ iff there is a unique {\it q.s.d.};
if $S=\infty$ then there is no {\it q.s.d.} or there are infinitely many ones. Then the equivalence between $(ii)$
and $(iv)$
is immediate. Also it is well known that
$$
\E_n(T_0)=\sum_{i\ge 1}\pi_i+\sum_{r=1}^{n-1}(\lambda_r \pi_r)^{-1}\sum_{i\ge r+1} \pi_i\,,
$$
see for example \cite[formula 7.10]{karlin}. Therefore (iii) and (iv) are equivalent.
Let us now
examine how this criterion is related to the nature of the boundary at $+\infty$. From the table
in \cite[section 8.1]{And} we have that $+\infty$ is an entrance boundary iff $A=\infty$, $\E_1(T_0)<\infty$
and $S<\infty$. Finally, $S<\infty$ implies that $\E_1(T_0)<\infty$ and this ensures $\P_1(T_0<\infty)=1$, that is
$A=\infty$. This shows the result.
\end{proof}

\section{\bf {Biological models}}

\subsection{Population dynamics and quasi-stationary distributions}

Our aim is to model the dynamics of an isolated population by a diffusion $Z:=(Z_t;t\ge 0)$.
Since competition for limited resources impedes natural populations with no
immigration to grow indefinitely, they are all doomed  to become extinct at some finite time
$T_0$. However, $T_0$ can be large compared to human timescale and it is common that population
sizes fluctuate for large amount of time before extinction actually occurs.
The notion of quasi-stationarity captures this behavior \cite{Ferri,SE04}.

The  diffusions we consider arise as scaling limits of general
birth--death processes. More precisely, let  $(Z^N)_{N\in \NN}$ be
a sequence of continuous time birth--death processes $Z^N:=(Z^N_t;
t\ge 0)$, renormalized by the weight $N^{-1}$, hence taking values
in $N^{-1}\mathbb{N}$. Assume that their birth and death rates
from state $x$ are equal to $b_N(x)$ and $d_N(x)$, respectively,
and $b_N(0)=d_N(0)=0$, ensuring that the state $0$ is absorbing.
We also assume that for each $N$ and for some constant $B_N$,
$$
b_N(x)\leq (x+1)B_N,\ x\geq 0
$$
and that there exist a nonnegative
constant $\gamma$ and a function $h \in C^1([0,\infty))$,
$h(0)=0$, called the \emph{growth function}, such that
\begin{equation}
\label{CVbndn} \forall x\in (0,\infty):\;\;\; \lim_{N\to \infty} \frac{1}{N}(b_N(x)-d_N(x)) =
h(x)\quad ; \quad \lim_{N\to \infty} \frac{1}{2N^2}(b_N(x)+d_N(x)) = \gamma x.
\end{equation}

Important ecological examples include
\begin{itemize}

\item[(i)] The \emph{pure branching} case, where the individuals give birth and die independently,
so that one can take  $b_N(x)=(\gamma N+\lambda)Nx$ and $d_N(x)=(\gamma N+\mu)Nx$.
Writing $r:=\lambda-\mu$ for the \emph{Malthusian growth parameter} of the population, one gets $h(z) = rz$.

\item[(ii)] The \emph{logistic branching} case, where $b_N(x)=(\gamma N+\lambda)Nx$ and
$d_N(x)=(\gamma N+\mu)Nx +\frac{c}{N} N x(N x-1)$. The quadratic
term in the death rate describes the interaction between
individuals. The number of individuals is of order $N$, the
biomass of each individual is of order $N^{-1}$, and $c/N$ is the
interaction coefficient. The growth function is then $h(z) =
rz-cz^2$.

\item[(iii)] Dynamics featuring \emph{Allee effect}, that is, a positive
density-dependence for certain ranges of density, corresponding to
cooperation in natural populations. A classical type of growth function in that setting is
$h(z)=rz({z\over K_0}-1)(1-{ z\over K})$. Observe that in this last case, the individual growth rate
is no longer a monotone decreasing function of the population size.
\end{itemize}

Assuming further that $(Z_0^N)_{N\in \NN}$ converges as
$N\to\infty$ (we thus model the dynamics of a population whose
size is of order $N$), we may  prove, following  Lipow \cite{LIP}
or using the techniques of  Joffe-M\'etivier \cite{JM86}, that the
sequence $(Z^N)_{N\in \NN}$ converges weakly to a continuous limit
$Z$. The parameter $\gamma$ can be interpreted as a demographic
parameter describing the ecological timescale. There is a main
qualitative difference depending on whether $\gamma=0$ or not.

If $\gamma=0$, then the limit $Z$ is a deterministic solution to
the dynamical system $ \dot{Z}_t=h(Z_t)$. Since $h(0)=0$, the
state $0$ is always an equilibrium, but it can be unstable. For
example, in the logistic case $h(z)= rz -cz^2$ and it is easily
checked that when $r>0$, the previous dynamical system has two
equilibria, $0$ which is unstable, and $r/c$ (called
\emph{carrying capacity}) which is asymptotically stable. In the
Allee effect case, $0$ and $K$ are both stable equilibria, but
$K_0$ is an unstable equilibrium, which means the population size
has a threshold $K_0$ to growth, below which it cannot take over.

If $\gamma>0$, the sequence $(Z^N)_{N\in \NN}$ converges in law to
the process $Z$, solution to the following stochastic differential
equation
\begin{equation}
\label{eqn : generalized Feller}
 dZ_t=\sqrt{\gamma Z_t}dB_t+h(Z_t)dt.
\end{equation}
The acceleration of the ecological process has generated the white noise.
Note that $h'(0^+)$ is the mean \emph{per capita} growth
rate for \emph{small} populations. The fact that it is finite is mathematically convenient, and
biologically reasonable. Since $h(0)=0$,
 the population undergoes no immigration, so
that 0 is an absorbing state.
 One can easily check that when time goes to infinity, either $Z$ goes
to $\infty$ or is absorbed at $0$.

When $h\equiv 0$, we get the classical Feller diffusion, so we call \emph{generalized Feller
diffusions} the diffusions driven by (\ref{eqn : generalized Feller}). When $h$ is
linear, we get the general continuous-state branching process with continuous paths,
sometimes also called Feller diffusion by extension. When $h$ is concave
quadratic, we get the logistic Feller diffusion \cite{Eth, Lam05}.

\medskip

\begin{definition}\label{def : hh}
(HH) We say that $h$ satisfies the condition (HH) if
$$
(i)\lim_{x \to \infty}\frac{h(x)}{\sqrt{x}}=-\infty,\qquad (ii)\lim_{x\to \infty}\frac{x
h'(x)}{h(x)^2} =0.
$$
\end{definition}

In particular (HH) holds for any subcritical branching diffusion, and any logistic Feller
diffusion. Concerning assumption $(i)$, the fact that $h$ goes to $-\infty$ indicates strong
competition in large populations resulting in negative growth rates (as in the logistic case).
Assumption $(ii)$ is fulfilled for most classical biological models, and
it appears as a mere technical condition.

\medskip
Gathering all results of the present paper and applying them to our biological model yields the
following statement, which will be proved at the end of this section.

\begin{theorem}\label{thm : main result}
Let $Z$ be the solution of \eqref{eqn : generalized Feller}. We assume $h\in C^1([0,\infty)), \; h(0)=0$ and 
that $h$ satisfies assumption (HH).
Then, for all initial laws with bounded support, the law of $Z_t$ conditioned on $\{Z_t\neq 0\}$
converges exponentially fast to a probability measure $\nu$, called the Yaglom limit.
The law $\Q_x$ of the process $Z$ starting from $x$ and conditioned to be never extinct exists and
defines the so-called $Q$-process. This process converges, as $t\to \infty$, in distribution, to
its unique invariant probability measure. This probability measure is absolutely continuous w.r.t.
$\nu$ with a nondecreasing Radon--Nikodym derivative.

In addition, if  the following integrability condition is satisfied
$$\int_1^{\infty}\frac{dx}{-h(x)}<\infty ,$$
then $Z$ comes down from infinity and the convergence of the conditional one-dimensional distributions
holds for all initial laws. In particular, the Yaglom limit $\nu$ is then the
unique quasi-stationary distribution.
\end{theorem}

\begin{proof}
For $Z$ solution to \eqref{eqn : generalized Feller}, recall that
$X_t = 2\sqrt{Z_t/\gamma}$ satisfies the SDE $dX_t =
dB_t -q(X_t) dt$ with
$$
q(x) = \frac{1}{2x}-\frac{2h(\gamma x^2 /4)}{\gamma x} \qquad x>0.
$$
In particular we have  $q'(x)=-\frac{1}{2x^2}+\frac{2h(\gamma x^2 /4)}{\gamma
x^2}-h'(\gamma x^2/4)$ and
$$
q^2(x)-q'(x)=\frac{3}{4 x^2} + h(\gamma x^2/4)\left(
\frac{4}{\gamma^2 x^2}h(\gamma x^2/4) -\frac{4}{\gamma
x^2}\right)+h'(\gamma x^2 /4).
$$
Under assumption (HH) we have the following behaviors at $0$ and $\infty$: 
$q(x)\mathop{\sim}\limits_{x\downarrow 0}1/2x$, as well as
$$
q^2(x)-q'(x)\mathop{\sim}\limits_{x\downarrow 0}\frac{3}{4x^2}
\qquad\mbox{ and }\qquad
(q^2-q')(2\sqrt{x/\gamma})\mathop{\sim}\limits_{x\to\infty}\frac{h(x)^2}{x}\left(\frac{1}{\gamma}+\frac{xh'(x)}{h(x)^2}\right).
$$
Then, it is direct to check that hypothesis (H2) holds
$$
\lim_{x \to \infty} q^2(x) - q'(x)=\infty\,
\textrm{ , } \, C:=-\inf_{x \in (0,\infty)} q^2(x)-q'(x) <\infty.
$$
We recall that 
$$Q(x)=\int_1^x 2q(y) dy, \; \Lambda(x)=\int_1^x e^{Q(y)} dy \hbox{ and } 
\kappa(x) = \int_1^x e^{Q(y)}\left(\int_1^y e^{-Q(z)}dz
\right) dy. 
$$
Straightforward calculations show that 
$$
\lim_{x\to \infty}\frac{Q(x)}{x}=\infty \qquad\mbox{ and
} \qquad A:=\lim_{x\rightarrow 0^+} \big(Q(x)-\log(x) \big)
\in(-\infty,\infty).
$$
In particular, $\Lambda(\infty) = \infty$ and the integrand in
the definition of $\kappa$ is equivalent to $y\log (y)$ which
ensures $\kappa(0^+)<\infty$. Thus $X$, and consequently $Z$,
is absorbed at 0 with probability 1, that is hypothesis (H1) holds.

We now continue with (H3) which is
$$
\int_0^1 \, \frac{1}{q^2(y) - q'(y) +
C + 2} \, e^{-Q(y)}dy < \infty.
$$
This hypothesis holds because near 0 the integrand is of the order
$$
\frac{1}{\frac{3}{4y^2}} e^{-Q(y)}\sim \frac{4e^{-A}}{3} y.
$$ 
For the first part of the Theorem it remains only to show that (H4) holds
$$
\int_{1}^{\infty}e^{-Q(x)}dx<\infty\quad\text{and}\quad
\int_{0}^{1} x\;e^{-Q(x)/2}dx<\infty.
$$
The first integral is finite because $Q$ grows at least linearly near $\infty$ and the second one
is finite because the integrand is of order $1/\sqrt{x}$ for $x$ near $0$.

Hence we can apply Theorem \ref{thm : Yaglom}, Proposition \ref{ratecv}, and
Corollaries \ref{corq} and \ref{corjavaitort} to finish with the proof of the 
first part of the Theorem.

\bigskip

For the last part of the Theorem we need to show that $X$ comes down from infinity
which is equivalent to  (H5). Thanks to Remark \ref{rempatrick}, there is a
simple sufficient condition for this hypothesis to hold, which has three
components. The first one 
$$
q(x)\geq q_0 > 0\; \hbox{ for all } x \geq x_0
$$
follows from (HH)$(i)$. The second one 
$$
\limsup_{x\to \infty} q'(x)/2 q^2(x)<1
$$
is equivalent to
$$
\limsup_{x\to \infty}-\frac{x h'(x)}{h(x)^2}
<\frac2\gamma,
$$
which clearly follows from (HH)$(ii)$.
Finally the third one
$$
\int_{x_0}^{\infty} \,
\frac{1}{q(x)} \, dx < \infty,
$$
thanks to (HH)$(i)$, is equivalent to
$$
\int^\infty \frac{-\gamma x}{2h(\gamma x^2/4)} \, dx=\int^\infty \frac{1}{-h(z)} \, dz<\infty.
$$
This is exactly the extra assumption made in the Theorem and the result is proven.
\end{proof}

\subsection{The growth function and conditioning}

\quad Referring to the previous construction of the generalized
Feller diffusion \eqref{eqn : generalized Feller}, we saw why
$h(z)$ could be viewed as the expected growth rate of a population
of size $z$ and $h(z)/z$ as the mean \emph{per capita} growth
rate. Indeed, $h(z)$ informs of the resulting action of density
upon the growth of the population, and $h(z)/z$ indicates the
resulting action of density upon each individual. In the range of
densities $z$ where $h(z)/z$ increases with $z$, the most
important interactions are of the \emph{cooperative} type, one
speaks of \emph{positive density-dependence}. On the contrary,
when $h(z)/z$ decreases with $z$, the interactions are of the
\emph{competitive} type, and density-dependence is said to be
\emph{negative}. In many cases, such as the logistic one, the
limitation of resources forces harsh competition in large
populations, so that, as $z\to \infty$,  $h(z)/z$ is negative and
decreasing. In particular $h(z)$ goes to $-\infty$. The shape of
$h$ at infinity determines the long time
behavior of the diffusion $Z$.\\

Actually, if $h$ goes to infinity at infinity, such as in the pure
branching process case (where $h$ is linear), Theorem \ref{thm :
main result} still holds if (HH)$(i)$ is replaced with the more
general condition $\lim_{x \to
\infty}\frac{h(x)}{\sqrt{x}}=\pm\infty$, \emph{provided} the
generalized Feller diffusion is further \emph{conditioned on
eventual extinction}. Indeed, the following statement ensures that
conditioning on extinction roughly amounts to \emph{replacing} $h$
\emph{with} $-h$.

\begin{proposition}
\label{thm : conditioning on ext} assume that $Z$ is given by
(\ref{eqn : generalized Feller}), where $h\in C^1([0,\infty))$,
$h(0)=0$, $\lim_{x \to \infty} \frac{h(x)}{\sqrt{x}}=\infty$. Define
$u(x):=\PP_x(\lim\limits_{t\to \infty}Z_t=0)$ and let $Y$ be the
diffusion $Z$ conditioned on eventual extinction. Then $Y$ is
the solution of the SDE, $Y_0=Z_0$
\begin{equation}
\label{sde*}
dY_t=\sqrt{\gamma Y_t}dB_t +\left(h(Y_t) + \gamma Y_t
\frac{u'(Y_t)}{u(Y_t)}\right)dt.
\end{equation}

If, in addition $h$
satisfies (HH)$(ii)$ then
$$
h(y) + \gamma y \frac{u'(y)}{u(y)}\sim_{y\to \infty} -h(y).
$$
\end{proposition}

\begin{proof}
Let $ J(x):=\int_0^x\frac{2h(z)}{\gamma z}dz$ which is
well-defined since $h\in C^1([0,\infty))$ with $h(0)=0$. We set
$$
v(x):=a\int_x^\infty e^{-J(z)}dz,
$$
with $a=(\int_0^\infty e^{-J(z)}dz)^{-1}$ (well-defined by the growth of $h$ near $\infty$).
Now we prove that $u=v$. It is easily
checked that $v$ is decreasing with $v(0)=1$, $v(\infty)=0$, and
that it satisfies the equation $\frac{\gamma}{2}xv''(x) +
h(x) v'(x)=0$  for all $x\geq 0$.

As a consequence, $(v(Z_t);t\ge 0)$ is a (bounded hence) uniformly
integrable martingale, so that
$$
v(x)=\EE_x(v(Z_t))\rightarrow v(0)\PP_x(\lim_{t\to
\infty}Z_t=0)+v(\infty)\PP_x(\lim_{t\to \infty}Z_t=\infty)=u(x),
$$
as $t\to\infty$, so that indeed $u=v$.

Using the strong Markov Property of $Z$ we obtain that for any Borel set
$A\subset (0,\infty)$ and $s\ge 0$
$$
\P_x(Y_s \in A)=\P_x(Z_s \in A \mid  T_0<\infty)=
\E_x\left(\frac{\P_{Z_s}(T_0<\infty)}{\P_x(T_0<\infty)},\, Z_s \in A\right)=
\E_x\left(\frac{u(Z_s)}{u(x)},\, Z_s \in A\right).
$$
Then for any measurable and bounded function $f$ we get
$$
\E_x(f(Y_s))=\E_x\left(f(Z_s) \frac{u(Z_s)}{u(x)}\right).
$$
Now if $f$ is $C^2$ and has compact support contained in $(0,\infty)$, we get
from It\^o's formula that $uf$ is in the domain of $L^Z$, the generator of $Z$,
and then $f$ is in the domain of the generator $L^Y$ of $Y$ and moreover
$$
L^Y(f)(x)=\frac{1}{u(x)}{L^Z(uf)(x)}=\frac{\gamma}{2}xf''(x)+\left(h(x)+\gamma
x\frac{u'(x)}{u(x)}\right)
$$
Then, since $h$ is locally Lipschitz we obtain that the law of $Y$ is the
unique solution to the SDE (\ref{sde*}).

\medskip

Let us show the last part of the proposition. Notice
that $J$ is strictly increasing after
some $x_0$, so we consider its inverse $\phi$ on $[J(x_0),
\infty)$. Next observe that for $x>x_0$,
$$
-\frac{u}{u'}(x) = e^{J(x)}\int_{x}^\infty e^{-J(z)} dz =
e^{J(x)}\int_{J(x)}^\infty e^{-b} \phi'(b)db,
$$
with the change $b=J(z)$. As a consequence, we can write for
$y>J(x_0)$
\begin{equation}
\label{eqn : referee relou 1} -\frac{u}{u'}(\phi(y)) =
e^{y}\int_{y}^\infty e^{-b} \phi'(b)db=\int_{0}^\infty e^{-b}
\phi'(y+b)db.
\end{equation}
Because $h$ tends to $\infty$, $J(x)\ge (1+\varepsilon)\log(x)$
for $x$ sufficiently large, so that $\phi(y)\le
\exp(y/(1+\varepsilon))$, and $\phi(y)\exp(-y)$ vanishes as $y\to
\infty$. Now, since
$$
\phi'(y)=\frac{\gamma \phi(y)}{2h(\phi(y))}=o(\phi(y)),
$$
$\phi'(y)\exp(-y)$ also vanishes. Since $h$ is differentiable, $J$
is twice differentiable, and so is $\phi$, so performing an
integration by parts yields
\begin{equation}
\label{eqn : referee relou 2} \phi'(y)=\int_{0}^\infty e^{-b}
\phi'(y+b)db\,-\,\int_{0}^\infty e^{-b} \phi''(y+b)db.
\end{equation}
Since $\phi'(J(x)) = 1/J'(x)$, we get $J'(x)\phi''(J(x)) = (1/J'(x))'$,
so by the technical assumption (HH)$(ii)$,
$$
\phi''(J(x)) = \phi'(J(x)) \left(\frac{1}{J'(x)}\right)'=
\frac{\gamma}{2}\phi'(J(x))\left(
\frac{1}{h(x)}-\frac{xh'(x)}{h(x)^2}\right)=o\left(\phi'(
J(x))\right),
$$
as $x\rightarrow\infty$. Then, as $y\rightarrow\infty$ we have $\phi''(y)=
o(\phi'(y))$. This shows, thanks to \eqref{eqn :
referee relou 2}, that
$$
\int_{0}^\infty e^{-b} \phi'(y+b)db \sim_{y\to \infty}\phi'(y)
$$
which entails, thanks to \eqref{eqn : referee relou 1}, that
$$
-\frac{u}{u'}(\phi(y)) \sim_{y\to \infty} \phi'(y).
$$
This is equivalent to
$$
\gamma x\frac{u'}{u}(x) \sim_{x\to \infty} -\gamma x J'(x)=
-2h(x),
$$
which ends the proof.
\end{proof}

Let us examine the case of the Feller diffusion (pure branching
process), where $h(z)=rz$. First, it is known (see e.g.
\cite[Chapter 2]{tojuanagua}) that when $r>0$, the supercritical
Feller diffusion $Z$ conditioned on extinction is \emph{exactly}
the subcritical Feller diffusion with $h(z)=-rz$. The previous
statement can thus be seen as an extension of this duality to more
general population diffusion processes.

Second, in the critical case ($r=0$), our present results do not
apply. Actually, the (critical) Feller diffusion has no {\it
q.s.d.} \cite{Lam06}. Third, in  the subcritical case ($r<0$), our
results do apply, so there is a  Yaglom limit and a $Q$-process,
but in contrast to the case when $1/h$ is integrable at $\infty$,
it is shown in \cite{Lam06} that  subcritical Feller diffusions
have infinitely many {\it q.s.d.}.


\appendix

\section{\bf {Proof of Lemma \ref{estimkernel}}}

We first prove the second  bound. For any nonnegative  and continuous function $f$  with support
in $\R^{+}$ we have from hypothesis (H2)
$$
\int \tilde p_{1}(x,u) \;f(u)du=
 \E^{\W_x} \left[ f(\omega(1)) \,
\BBone_{1<T_0}(\omega) \, \exp\left(-  \, \frac 12 \, \int_0^1 \, (q^2 - q')(\omega_s)
ds\right)\right]
$$
$$
\le e^{C/2} \E^{\W_x} \left[ f(\omega(1)) \, \BBone_{1<T_0}(\omega)\right]\;.
$$
The  estimate (\ref{estim2}) follows by letting $f(z)dz$ tend to the Dirac
measure at $y$ with $K_3=e^{C/2}$, that is
$$
\tilde p_{1}(x,y)\le K_3\; p_{1}^{D}(x,y).
$$
Here $p_{1}^{D}(x,y)=\frac{1}{\sqrt{2\pi}}\left(e^{-(x-y)^2/2}-e^{-(x+y)^2/2}\right)$
(see for example \cite[page 97]{karatzas}).

\medskip

Let us now prove the
upper bound in  (\ref{estim1}). Let $B_{1}$ be the function defined by
$$
B_{1}(z):=\inf_{u\ge z} \big(q^{2}(u)-q'(u)\big)\;.
$$
We have
$$
\int \tilde p_{1}(x,y) \;f(y)dy= \E^{\W_x} \left[ f(\omega(1)) \, \BBone_{1<T_0}
\BBone_{1<T_{x/3}}\, \exp\left(-  \, \frac 12 \, \int_0^1 \, (q^2 - q')(\omega_s) ds\right)\right]
$$
$$
+\E^{\W_x} \left[ f(\omega(1)) \, \BBone_{1<T_0} \, \BBone_{1\ge T_{x/3}}\exp \left(-  \, \frac 12
\, \int_0^1 \, (q^2 - q')(\omega_s) ds\right)\right]\;.
$$
For the first expectation we have
$$
\E^{\W_x} \left[ f(\omega(1)) \, \BBone_{1<T_0} \BBone_{1<T_{x/3}}\, \exp\left(-  \, \frac 12 \,
\int_0^1 \, (q^2 - q')(\omega_s) ds\right)\right]
$$
$$
\le e^{-B_{1}(x/3)/2}\; \E^{\W_x} \left[ f(\omega(1)) \, \BBone_{1<T_0}\right]\;.
$$
For the second expectation, we obtain
$$
\E^{\W_x} \left[ f(\omega(1)) \, \BBone_{1<T_0} \, \BBone_{1\ge T_{x/3}}\exp \left(-  \, \frac 12
\, \int_0^1 \, (q^2 - q')(\omega_s) ds\right)\right]
$$
$$
\le e^{C/2} \;\E^{\W_x} \left[ f(\omega(1)) \, \BBone_{1<T_0} \BBone_{1\ge T_{x/3}}\right]
$$
$$
=
 e^{C/2} \left(\E^{\W_x} \left[ f(\omega(1)) \,
\BBone_{1<T_0}\right] - \;\E^{\W_x} \left[ f(\omega(1)) \, \BBone_{1< T_{x/3}}\right]\right)\;.
$$
Using a limiting argument as above and the invariance by translation of the law of the Brownian
motion, and firstly assuming that $y/2<x<2y$, we obtain
$$
 \tilde p_{1}(x,y)\le  e^{-B_{1}(x/3)/2} p^{D}_{1}(x,y)+
e^{C/2}\left(p^{D}_{1}(x,y)-p^{D}_{1}(2x/3,y-x/3)\right)\;,
$$

$$
p^{D}_{1}(x,y)-p^{D}_{1}(2x/3,y-x/3)=\frac{1}{\sqrt{2\pi}}
\left(e^{-(y+x/3)^{2}/2}-e^{-(x+y)^{2}/2}\right) \le
\frac{1}{\sqrt{2\pi}}\;e^{-\max\{x,y\}^{2}/18}\;.
$$
Since the function $B_{1}$ is non-decreasing, we  get for $y/2<x<2y$
$$
 \tilde p_{1}(x,y)\le   \frac{1}{\sqrt{2\pi}}\left(
 e^{-B_{1}(\max\{x,y\}/6)/2}+e^{-\max\{x,y\}^{2}/18}\right)\;.
$$
If $x/y\notin]1/2,2[$, we get from the estimate (\ref{estim2})
$$
\tilde p_{1}(x,y)\le  \frac{K_3}{\sqrt{2\pi}} \;e^{-(y-x)^{2}/2}\le \frac{K_3}{\sqrt{2\pi}}
e^{-\max\{x,y\}^{2}/8}\;.
$$
We now define the function $B$ by
$$
B(z):=\log\left(\frac{K_3\vee 1}{\sqrt{2\pi}}\right)+\min\big\{B_{1}(z/6)/4\;,\;z^{2}/36\big\}\;.
$$
It follows from hypothesis (H2) that $\lim_{z\to\infty}B(z)=\infty$. Combining the previous
estimates we get for any $x$ and $y$ in $\R^{+}$
$$
\tilde p_{1}(x,y)\le e^{-2B(\max\{x,y\})}\;.
$$
The upper estimate (\ref{estim1}) follows by taking the geometric average of this result and
(\ref{estim2}). We now prove that $\tilde p_{1}(x,y)>0$. For this purpose, let $a=\min\{x,y\}/2$
and $b=2\max\{x,y\}$. We have as above for every  nonnegative continuous function $f$ with support
in $\R^{+}$
$$
\int \tilde p_{1}(x,y) \;f(y)dy\ge
 \E^{\W_x} \left[ f(\omega(1)) \,
\BBone_{1<T_{[a,b]}} \, \exp\left(-  \, \frac 12 \, \int_0^1 \, (q^2 - q')(\omega_s)
ds\right)\right]
$$
where we denote $T_{[a,b]}$ the exit time from the interval $[a,b]$. Let
$$
R_{a,b}=\sup_{x\in[a,b]}( q^2(x) - q'(x))\;,
$$
this quantity is finite since $q\in C^{1}((0,\infty))$. We obtain immediately
$$
\int \tilde p_{1}(x,y) \;f(y)dy\ge e^{-\,R_{a,b}/2}\int p^{[a,b]}_{1}(x,y) \;f(y)dy
$$
where we denote $p^{[a,b]}_{t}$ the heat kernel with Dirichlet conditions in $[a,b]$. The
result follows from a limiting argument as above since $p^{[a,b]}_{1}(x,y)>0$. $\square$

\bigskip

\noindent {\bf Acknowledgements} S. Mart\'inez and J. San Mart\'in
thank the support from Nucleus Millennium P04-069-F. S.
M\'el\'eard also thanks the support from ECOS-Conicyt C05E02. The
authors thank the referee for a very thorough and
careful reading of the paper, as well as many helpful comments and
suggestions.

\bigskip

\bibliographystyle{plain}

\end{document}